\documentclass[11pt]{article}

\usepackage{fullpage}
\usepackage[T1]{fontenc}
\usepackage{amsmath,amssymb,amsfonts,amsthm,mathtools}
\usepackage{mathrsfs}
\usepackage{graphicx}
\usepackage{bm}
\usepackage{enumitem}
\usepackage{tikz}
\usetikzlibrary{matrix}
\usepackage{diagbox}
\usepackage{algorithm}
\usepackage{algpseudocode}
\usepackage{aliascnt}
\usepackage[hidelinks]{hyperref}
\usepackage[capitalise]{cleveref}
\graphicspath{{./}{plots/}}

\numberwithin{equation}{section}
\numberwithin{algorithm}{section}

\theoremstyle{plain}
\newtheorem{theorem}{Theorem}[section]
\newaliascnt{lemma}{theorem}
\newtheorem{lemma}[lemma]{Lemma}
\aliascntresetthe{lemma}
\newaliascnt{proposition}{theorem}
\newtheorem{proposition}[proposition]{Proposition}
\aliascntresetthe{proposition}
\newaliascnt{corollary}{theorem}

\aliascntresetthe{corollary}
\newaliascnt{definition}{theorem}
\newtheorem{definition}[definition]{Definition}
\aliascntresetthe{definition}
\newaliascnt{assumption}{theorem}

\aliascntresetthe{assumption}
\theoremstyle{remark}
\newaliascnt{remark}{theorem}
\newtheorem{remark}[remark]{Remark}
\aliascntresetthe{remark}
\crefname{equation}{}{}
\Crefname{lemma}{Lemma}{Lemmas}
\Crefname{proposition}{Proposition}{Propositions}

\newcommand{\lc}{\left(}
\newcommand{\rc}{\right)}

\newcommand{\EE}{\mathbb{E}}

\newcommand{\N}{\mathbb{N}}

\newcommand{\mC}{\mathcal C}

\def\su{{\mathsf{u}}}
\def\sv{{\mathsf{v}}}

\def\bjm{{\bar{\jmath}}}

\newcommand{\ux}{{\underline x}}

\newcommand{\beq}{\begin{equation}}
\newcommand{\eeq}{\end{equation}}

\newcommand{\g}{\gamma}

\newcommand{\p}{\psi}

\newcommand{\D}{\Delta}
\newcommand{\Dx}{{\Delta x}}

\title{Comparison Principle and Semi-Lagrangian Methods for a Hamilton--Jacobi--Bellman System with Recursive Utility}
\author{Qing Tang\thanks{School of Mathematics and Physics, China University of Geosciences (Wuhan), Wuhan, China ({tangqingthomas@gmail.com}).}}
\date{\today}

\begin{document}
\maketitle

\begin{abstract}
We consider a system of coupled first-order Hamilton--Jacobi--Bellman equations arising in continuous-time heterogeneous agent models with recursive utility, a class of models that has attracted increasing attention in macroeconomics. The first part of the paper establishes a comparison principle for the system when preferences exhibit early resolution of uncertainty. A key ingredient in the proof is an ordinally equivalent transformation of the recursive utility representation. In the second part, we develop semi-Lagrangian schemes and iterative algorithms for solving the associated nonlinear discrete systems. For preferences exhibiting late resolution of uncertainty, we establish a first-order error estimate for the semi-Lagrangian approximation. Numerical experiments demonstrate that the proposed methods yield results consistent with those obtained using the previously developed finite difference schemes.


\end{abstract}

\section{Introduction}
In this paper, we consider a system of coupled first order Hamilton--Jacobi--Bellman equations with recursive utility, recently analyzed in
\cite{achdou2026continuous,achdou2026finite}.  Such a system of equations arises from the Aiyagari--Bewley--Huggett type heterogeneous agent 
(\textit{HA} for brevity) models in macroeconomic theory \cite{ljungqvist2018recursive}, extended to continuous time in
\cite{achdou2022income}. \textit{HA} models provide a standard framework for studying consumption--saving decisions, wealth distributions, and equilibrium asset prices in incomplete markets. By incorporating idiosyncratic income risk and borrowing constraints, these models capture the interaction between individual saving behavior and the distribution of wealth in general equilibrium. Their numerical solution involves 
Hamilton--Jacobi--Bellman (HJB) equations coupled with equations
describing the invariant distribution of agents and the equilibrium
conditions for asset prices. \par
The HJB equation for characterizing the agent's consumption--saving decision involves the Epstein--Zin recursive utility, which is specified through the aggregator
(cf. \cite{wang2016optimal})
\begin{equation}
\label{EZ}
f(c,v)=\frac{\rho}{1-\p^{-1}}
\frac{c^{1-\p^{-1}}-((1-\g)v)^\theta}
{((1-\g)v)^{\theta-1}},
\qquad
\theta=\frac{1-\p^{-1}}{1-\g}.
\end{equation}
In \cref{EZ}, $c$ denotes consumption, $v$ is the value function, $\rho$ is the
subjective discount rate, $\psi$ is the elasticity of intertemporal
substitution (EIS), and $\gamma$ is the coefficient of relative risk
aversion. We focus on the parameter region $0<\psi<1$ and $\gamma>1$.
The time-additive constant relative risk aversion (CRRA) specification is recovered when
$\gamma=\psi^{-1}$. The parameter $\theta$ determines the preference for the
timing of the resolution of uncertainty: early resolution is preferred when
$\theta<1$, while late resolution is preferred when $\theta>1$
\cite{EZ1989}.\par 
The system of HJB equations with recursive utility and a given interest rate $r$ is 
\begin{equation}\label{Eq: HJB}
0 = \max_{c\ge 0} \left\{ f(c,v_j) +(rx+y_j - c)Dv_j(x) \right\}+\lambda_j(v_{\bar \jmath}(x)-v_j(x)),\quad j\in \{1,2\},\,\,\bar\jmath=3-j.
\end{equation}
The equations are coupled through the  switching terms $\lambda_j(v_{\bar\jmath} - v_j)$, encoding the transitions between income states with Poisson rates $\lambda_j$. The state constraint $x \geq \ux$ (a borrowing limit) is imposed to prevent negative savings at $\ux$: $r\ux+y_j - c_j(\ux)\geq 0$.  
 \par
 As a first theoretical contribution, we prove a comparison principle for solutions of \cref{Eq: HJB} in the case of preference for early resolution of uncertainty, i.e. $0<\theta<1$. The uniqueness of solutions to this type of problem in the discrete-time setting has been analyzed in \cite{peng2026dynamic,ren2021dynamic,Sargent2025}. In the continuous-time setting, the strong comparison principle (between semicontinuous sub- and supersolutions) in the case $\theta= 1$ has been established in \cite{SL_HA} and then extended to $\theta>1$ in \cite{achdou2026continuous}. The main difficulty in the case $0<\theta<1$ is that the nonlinear part of the aggregator $f(c,v)$ in \cref{Eq: HJB} becomes an increasing function for the $v-$argument. The existence of extremal solutions has been addressed in \cite{achdou2026finite} but the uniqueness of solution was not established. In this paper we exploit an ordinal transformation of the recursive utility problem to obtain a more tractable equivalent representation. Because the transformation is invertible, comparison principle in the transformed representation translates into that of the original problem. This yields a self-contained proof of uniqueness to the HJB equation in the parameter regime $0<\theta<1$. Similar changes of variables have been used in \cite{alvarez1996viscosity,duffie1992pde} for second order HJB equation involving recursive utility. 
\par 
 We next address the numerical solution of such HJB systems and related \textit{HA} models. The finite difference methods for solving \cref{Eq: HJB} have been considered recently in \cite{achdou2026finite}.
Building on the semi-Lagrangian framework for state-constrained heterogeneous agent models with time-additive CRRA utility developed in \cite{SL_HA}, we extend the HJB discretization to Epstein--Zin recursive utility, whose nonlinear dependence on the value function requires new analytical arguments and iterative algorithms. SL methods discretize the dynamic programming principle by following the controlled dynamics over a discrete time interval and evaluating the value function at the resulting off-grid points through interpolation. We refer to \cite{bardi1997optimal,falcone2013semi} for the theoretical foundations, to \cite{briani2012approximation} for approximation of coupled HJB systems, to \cite{carlini2020semi} for convergence rate analysis and to \cite{carlini2014fully,ashrafyan2025fully} for SL schemes for Mean Field Games. Concerning iterative algorithms, \cite{carlini2026semi} studied a learning value iteration algorithm and its acceleration by policy iteration for a semi-Lagrangian discretization of Mean Field Games. In the abstract dynamic programming framework, \cite{Sargent2025,peng2026dynamic} developed monotonicity-based and order-theoretic approaches to establish convergence of iterative algorithms, for solving problems with recursive preferences.\par
We consider the underlying stochastic control problem associated with the Epstein--Zin utility, cf.
\cite{duffie1992stochastic,duffie1992pde,wang2016optimal}:
\begin{equation}\label{control_problem}
	v(x,y)=\sup_{c_\tau}\EE \left[\int_t^{\infty} f(c_\tau,v_\tau)\, d\tau \big\vert \boldsymbol{x}_t=x,\,\boldsymbol{y}_t=y\right],\, \text{s.t.}\,
      \left\{
        \begin{array}[c]{rcl}
          d\boldsymbol{x}_\tau &=& (r \boldsymbol{x}_\tau+\boldsymbol{y}_{\tau} - c_\tau)d\tau,\,\,\tau>t,\\
          \boldsymbol{x}_\tau  &\geq& \underline{x}.    
        \end{array}
      \right.  
	\end{equation}
	Here, $\boldsymbol{x}_\tau$ and $\boldsymbol{y}_{\tau}$ respectively stand for the agent's wealth and labor income at time $\tau$. The control variable is the consumption $c_\tau$ and the flow $f$ depends also on the continuation value $v_\tau=v(\boldsymbol{x}_\tau,\boldsymbol{y}_{\tau})$. It is assumed that $\boldsymbol{y}_{\tau}$ is a two-state continuous-time Markov chain with Poisson switching rates,. 
We observe that the same ordinal equivalent transformation for the HJB equation can also apply to the SL scheme, which allows us to establish a discrete comparison principle in the case $0<\theta<1$. Building on the Howard policy iteration in \cite{SL_HA}, we develop a Howard--Newton algorithm in which the nonlinear policy evaluation arising from recursive utility is performed using Newton's method. We prove convergence of the outer policy iteration in both the late- and early-resolution preference regimes. As an alternative, in the
early-resolution case $0<\theta<1$, we formulate the discrete problem as a
fixed point of a monotone and invariant map and develop a corresponding
value iteration algorithm. \par
This paper is organized as follows. In \cref{Sec: prelim}, we introduce the relevant background and previous results for the HJB equation, building on \cite{achdou2026continuous, achdou2026finite}. In \cref{Sec: comparison}, we prove the comparison principle and uniqueness of solution to the HJB equation in the case $0<\theta<1$ (a result which is not available in \cite{achdou2026continuous, achdou2026finite}), based on an ordinal equivalent transformation. In \cref{Sec: SL 1}, we propose an SL scheme. We develop a Howard--Newton algorithm and prove a discrete comparison principle for the solution to the scheme. In the early-resolution case, we propose additionally a monotone value iteration algorithm. For the late-resolution case, the SL scheme has a first-order convergence rate to the unique solution of the HJB equation. In \cref{Sec: MFG}, we give the solution methods for the Fokker--Planck--Kolmogorov (FPK) equations and a general equilibrium Mean Field Game model. Finally, numerical experiments in \cref{Sec: numerical} demonstrate the effectiveness of the proposed methods and confirm their agreement with the results obtained using the finite difference methods from \cite{achdou2026finite}. 

\section{Preliminaries}\label{Sec: prelim}
The assumptions that follow will be made in the whole paper:
\beq\label{asp1}
\g>1,\,\,0<\p<1,\quad \rho> r>0,\quad y_2>y_1>0,\quad\ux\leq 0,\quad  \rho \ux+y_1>0.
\eeq
It has been proved in \cite{achdou2026continuous} that the equilibrium interest rate $r^*$ is smaller than $\rho$, if there exists a solution to \cref{MFG} below. This justifies the bounds on $r$ in \cref{asp1}. With the constant $b$ defined in \cref{tab1} below, \cref{asp1} yields $r<b<\rho$, $rx+y_2>rx+y_1>0$ and $b(x+y_1/r)>0$. Therefore, the functions in \cref{Eq: sub super} below are well defined. \par  
We observe that the Epstein-Zin utility can be decomposed as, using $\rho/\theta=\zeta$,
	\beq\label{EZ flow}
	f(c,v_j)=\mathcal{F}(c,v_j)-\zeta v_j,\qquad \mathcal{F}(c,v_j)=\frac{\rho}{1-\p^{-1}}\frac{c^{1-\p^{-1}}}{((1-\g)v_j)^{\theta-1}}.
	\eeq
	 We use \cref{EZ flow} to rewrite the HJB equation \cref{Eq: HJB} as  
	\beq\label{HJB}
	\begin{aligned}
	\frac{\rho}{\theta} v_j(x)=	H(x,y_j,v_j,Dv_j)+\lambda_j(v_{\bar \jmath}(x)-v_j(x)),
	\end{aligned}
	\eeq
	where we use the notation from \cref{EZ flow} and
\beq\label{Eq: H from F}
 H(x,y_j,v_j,p):=\max_{c\ge 0} \left\{ \mathcal{F}(c,v_j) +(rx +y_j- c)p\right\}.
\eeq
From the first order necessary optimality condition, the Hamiltonian in \Cref{Eq: H from F} is
	\begin{equation}\label{def H}
  H(x,y,v,p)=\left\{
      \begin{aligned}
     \quad  (rx +y)p+\frac{\rho^{\p}}{\p-1} p^{1-\p}((1-\g)v)^{\frac{1-\g \p}{1-\g}},\quad \quad & \text{if }\quad p\geq 0,\\
       +\infty, \quad \quad  & \text{if }\quad p< 0.
     \end{aligned}
   \right.
\end{equation}
We observe that $H(x,y,v,p)$ defined in \Cref{def H} is strictly convex in $p$ for fixed $(x,y,v)$ with $p>0$.
The optimal consumption (away from the borrowing limit) is given by the first order necessary optimality condition, whenever $Dv_j >0$, 
	\beq\label{c optimal}
	c_j=\mathop{\arg\max}\limits_{c\geq 0} \left\{ \mathcal{F}(c,v_j) -cDv_j \right\}=\rho^{\p}(Dv_j)^{-\p}((1-\g)v_j)^{\frac{1-\g \p}{1-\g}}.
	\eeq

We summarize the notation in \Cref{tab1}. 
\begin{table}[h!]
    \centering
    \caption{Symbols}
    \label{tab1}
    \begin{tabular}{l c c}
        \hline
        Subjective discount rate & $\rho>0$\\
        Risk aversion  & $\gamma>1$  \\
        Elasticity of intertemporal substitution (EIS) &$0<\psi<1$\\
        A parameter arising in Epstein-Zin utility &$\theta:=\frac{1-\p^{-1}}{1-\g}>0$\\
        A convenient parameter &$\zeta:= \rho/\theta$\\
        A convenient parameter &$b:= \rho\left[\frac{r+\p(\rho-r)}{\rho}\right]^{\frac{1}{1-\p}}$\\
       \hline
    \end{tabular}
\end{table}\par
For any $c>0$, if $\theta > 1$ (resp. $0<\theta<1$) then $\mathcal{F}(c,v)$ is decreasing (resp. increasing) in $v$, i.e. 
\beq\label{F monotone}
\mathcal{F}_v(c,v)< 0\quad (resp. \,\,\,\mathcal{F}_v(c,v)>0).
\eeq
If $\theta > 1$, then $f(c,v)$ and $\mathcal{F}(c,v)$ are jointly concave in $(c,v)\in (0,+\infty)\times (-\infty,0)$.
\begin{definition}\hfill\label{def vis sol}
	\begin{itemize}
		\item[1.] An upper semicontinuous (u.s.c.) function $v=(v_1,v_2)$ is said to be a viscosity subsolution of \eqref{HJB} at $x$, if whenever $\varphi$ is a smooth function and $v_j-\varphi$ has a local maximum at $x$, then
	\beq\label{Eq: sub}	
	\frac{\rho}{\theta} v_j(x)\le H(x,y_j,v_j(x),D\varphi(x) )+\lambda_j (v_{\bar \jmath}(x)-v_j(x)).
	\eeq	
		\item[2.] A lower semicontinuous (l.s.c.) function $v=(v_1,v_2)$ is said to be a viscosity supersolution of \eqref{HJB} at $x$, if whenever $\varphi$ is a smooth function and $v_j-\varphi$ has a local  minimum at $x$, then
	\beq\label{Eq: super}
	\frac{\rho}{\theta} v_j(x)\ge H(x,y_j,v_j(x),D\varphi(x) )+\lambda_j (v_{\bar \jmath}(x)-v_j(x)).
	\eeq	
	\item[3.] 	A continuous function $v$ is said to be a constrained viscosity solution to  system \eqref{HJB} if $v$ is a viscosity supersolution in $(\ux,\infty)$ and a viscosity subsolution in $[\ux,\infty)$.
	\end{itemize}
\end{definition}
The sub- and supersolutions proposed below will play an important role in what follows.  
\begin{proposition}\label{Prop:sub-super}
With the constant $b$ defined in \cref{tab1}, let $\bar{y}>y_2$ and consider the functions
\beq\label{Eq: sub super}
\begin{aligned}
\check{\su}_1(x)=\check{\su}_2(x)=\frac{(rx+y_1)^{1-\gamma}}{1-\gamma},\quad
 \check{\sv}_1(x)=\check{\sv}_2(x)= \frac{(b(x+\bar{y}/r))^{1-\gamma}}{1-\gamma}.
\end{aligned}
\eeq
 The pairs $(\check{\su}_1,\check{\su}_2)$ and $ (\check{\sv}_1,\check{\sv}_2)$ are respectively a subsolution of \eqref{HJB} in $[\ux,+\infty)$ and a supersolution of \eqref{HJB} in $(\ux,+\infty)$.
\end{proposition}
\begin{proof}
For the subsolution property, we refer to \cite[Proposition 3.5]{achdou2026continuous} for the case $\theta\geq 1$, noting that the same argument extends directly to the case $0<\theta<1$. We now prove that $(\check{\sv}_1,\check{\sv}_2)$ satisfies \cref{Eq: super} with $j=2$. Observe that $\check{\sv}_2\in C^1$, hence for all $x>\ux$, we know $D\check{\sv}_2(x)=D\varphi(x)$ if $\check{\sv}_2-\varphi$ has a local  minimum at $x$, 
$$
\begin{aligned}
&H(x,y_2,\check{\sv}_2(x),D\check{\sv}_2(x))+\lambda_2 (\check{\sv}_1(x)-\check{\sv}_2(x))\\
={}&\max_{c\ge 0} \left\{ \mathcal{F}\lc c,\check{\sv}_2(x)\rc +(rx +\bar{y}- c)D\check{\sv}_2(x)\right\}+(y_2-\bar{y})D\check{\sv}_2(x)
\leq \frac{\rho}{\theta}\check{\sv}_2(x).
\end{aligned}
$$
Similarly, we can derive \cref{Eq: super} with $j=1$. Therefore, $ (\check{\sv}_1,\check{\sv}_2)$ is a supersolution. 
\end{proof}
\begin{remark}
The sub- and supersolutions in \cref{Prop:sub-super} are similar to those constructed in \cite{achdou2026continuous} and \cite[Proposition 2.2]{achdou2026finite}. The only difference is that, in the present paper, we take $\bar{y}>y_2$ for the supersolution, whereas \cite[Proposition 2.2]{achdou2026finite} takes $\bar{y}=y_2$. This modification allows us to use analogous grid functions to construct suitable barriers for the discrete solution of the SL scheme, as will become clear in \cref{Prop: disc barriers SL} below.
\end{remark}
The strong comparison principle for viscosity sub- and supersolutions of \eqref{HJB} has been shown in \cite[Theorem 2.3]{SL_HA} with $\theta=1$ and \cite[Proposition 3.1]{achdou2026continuous} with $\theta>1$. Such results yield the uniqueness of the viscosity solution when $\theta \geq 1$. Reference \cite{achdou2026continuous} contains a complete study of continuous-time {\it{HA}} models in the case $\theta\geq 1$. The following proposition summarizes results from \cite[Proposition 3.5, Proposition 3.10, Proposition 3.11]{achdou2026continuous} .

\begin{proposition}\label{Prop: v sol}
Assume $\theta \geq 1$. There exists a unique viscosity solution $(v_1,v_2)$ which is $C^1$ and strictly concave. Moreover, $v_j\in W^{2,\infty}_{loc}(\ux,+\infty)$ and $v_2\geq v_1$. 
\end{proposition}
Next we give some results concerning the saving policies $s_j$. The following proposition summarizes results from \cite[Proposition 3.12, Corollary 3.2]{achdou2026continuous}.
\begin{proposition}\label{Prop: saving}
Assume $\theta \geq 1$. The optimal saving policy $s_1$ has the following properties: $s_1(x)< 0$ for all $x> \ux$ and $s_1(\ux)=0$. If furthermore
\beq\label{rho-r2}
\lc \rho-r\rc(r\ux+y_2)^{-\g}+\lambda_2\lc (r\ux+y_2)^{-\g}-(r\ux+y_1)^{-\g}\rc< 0,
\eeq
then $s_2(\ux)>0$. Moreover, $Dv_1(\ux)>Dv_2(\ux)$. 
\end{proposition}

The next result deals with the behavior of $s_2$ as $x\to +\infty$. It justifies solving the HJB equation numerically on a bounded domain $[\ux,\bar{x}]$. 
\begin{proposition}[\cite{achdou2026continuous}, Proposition 3.15]\label{Prop: s2 barx}
Assume $\theta \geq 1$ and $s_2(\ux)>0$. There exists $\widehat{x}>\ux$ such that $s_2(\widehat{x})=0$ and $s_2(x)<0$ for all $x> \widehat{x}$.
\end{proposition}
\begin{remark}
\cref{Prop: s2 barx} implies that if $\bar{x}>\widehat{x}$, then imposing an artificial state constraint at $\bar{x}$ while solving \cref{HJB} does not change the solution.
\end{remark}
It has been shown in \cite[Corollary 3.3]{achdou2026continuous} that, under the same assumptions as \cref{Prop: saving},  $\lim_{x\to \ux}D^2v_1(x)=-\infty$. The next result, from \cite[Appendix C]{achdou2026continuous}, is about the asymptotic analysis of $s_1$, $D^2v_1$ and $Dv_1$ near $x=\ux$. 
\begin{proposition}\label{Prop: Dv1 ux}
Assume $\theta \geq 1$.  We define 
\beq\label{Eq: def vark}
\varkappa:=\lc \frac{\rho}{\theta}-r\rc Dv_1(\ux)+\lambda_1(Dv_1(\ux)-Dv_2(\ux))-H_v(\ux,y_1,v_1(\ux),Dv_1(\ux))Dv_1(\ux).
\eeq
Near $x=\ux$, the functions $s_1$, $D^2v_1$ and $Dv_1$ have the following behavior: 
\beq\label{Eq: s1D2v1 ux}
s_1(x)D^2v_1(x)= \varkappa+o(1),
\eeq
\beq\label{Eq: -s1 bound}
s_1(x)=-\sqrt{\frac{2\varkappa \p}{\rho} \lc r\ux+y_1\rc^{1+\p^{-1}}((1-\g)v_1(\ux))^{\frac{\p^{-1}-\g}{\g-1}}(x-\ux)}+o(\sqrt{x-\ux}),
\eeq
\beq\label{Eq: D2v1 ux}
D^2v_1(x)= -\varkappa \lc\frac{2\varkappa \p}{\rho} \lc r\ux+y_1\rc^{1+\p^{-1}}((1-\g)v_1(\ux))^{\frac{\p^{-1}-\g}{\g-1}}(x-\ux)\rc^{-1/2}(1+o(1)),
\eeq
\beq\label{Eq: Dv1 ux}
Dv_1(x)= \rho\frac{((1-\g)v_1(\ux))^{\frac{\p^{-1}-\g}{1-\g}}}{ (r\ux+y_1)^{\p^{-1}}}-\sqrt{\frac{2\rho \varkappa(x-\ux)}{\p \lc r\ux+y_1\rc^{1+\p^{-1}}} ((1-\g)v_1(\ux))^{\frac{\p^{-1}-\g}{1-\g}}}+o(\sqrt{x-\ux}).
\eeq
\end{proposition}
Finally, we observe that the existence of a viscosity solution to \cref{Eq: HJB} has been proved in \cite[Proposition 2.15]{achdou2026finite} with a fixed point method.
\begin{proposition}\label{Prop: v sol early}
Assume $0<\theta <1$. There exists a viscosity solution $(v_1,v_2)$ to the system \cref{Eq: HJB}. 
\end{proposition}

\section{Comparison principle and uniqueness of solution in the case \texorpdfstring{$0<\theta<1$}{0 < theta < 1}}\label{Sec: comparison}
We perform the change of variable
\begin{equation}\label{Eq: v to w}
(1-\g)v_j = (-w_j)^{1/\theta}, \qquad j=1,2,
\end{equation}
equivalently
\begin{equation}\label{cov2}
w_j := -\big((1-\g)v_j\big)^{\theta}.
\end{equation}
Since $1-\g<0$, $\theta>0$, and $v_j<0$, we have $w_j<0$. Moreover $u\mapsto -((1-\g)u)^\theta$ is an increasing transformation of $u$, so $w_j$ inherits monotonicity in $x$ from $v_j$.\par
The following elementary function is important for the analysis of the HJB equation in the case $0<\theta<1$. For $a,b>0$ we define
\begin{equation}\label{Phi-def}
\Phi(a,b) := a^{1-1/\theta}\,b^{1/\theta}.
\end{equation}
This will play the role of the switching/coupling nonlinearity between the two income states.\par
\begin{lemma}\label{lem:Phi-convex}
Assume $0<\theta<1$, then $\Phi(a,b)=a^{1-1/\theta}b^{1/\theta}$ is jointly convex on $(0,\infty)^2$.
\end{lemma}
\begin{lemma}\label{lem:Phi-shift}
For $a,b>0$ and $\delta\ge0$,
\begin{equation}\label{lem:Phi-shift}
\Phi(a+\delta,\,b+\delta) - \Phi(a,b) \ \le\ \delta.
\end{equation}
\end{lemma}

\begin{proof}
We first observe that joint convexity together with positive homogeneity of degree $1$ (i.e.\ $\Phi(ta,tb)=t\Phi(a,b)$ for $t>0$), 
$
 \frac{1}{2}\Phi(a+\delta,\,b+\delta) =\Phi\Big(\frac{a+\delta}{2},\frac{b+\delta}{2}\Big) \ \le\ \frac{1}{2}\Phi(a,b)+\frac{1}{2}\Phi(\delta,\delta).
$
By \eqref{Phi-def} and $\big(1-\tfrac1\theta\big)+\tfrac1\theta=1$, $\Phi(\delta,\delta)=\delta^{1-1/\theta}\delta^{1/\theta}=\delta$, which yields \cref{lem:Phi-shift}.
\end{proof}
Since the transformation is smooth and strictly monotone increasing, we observe that if $(v_1,v_2)$ is a viscosity solution to \cref{Eq: HJB}, then $(w_1,w_2)$ is a viscosity solution to the following system of HJB equations: 
\begin{equation}\label{Eq: HJB-w}
\begin{aligned}
\rho\, w_j(x) ={}& \sup_{c\ge0}\Big\{ (rx+y_j-c)\,Dw_j(x) - \rho\, c^{1-\p^{-1}} \Big\}
\ -\ \theta\lambda_j\Big[w_j(x) + \Phi\big(-w_j(x),\,-w_{\bar\jmath}(x)\big)\Big]\\
={}&\mathcal{H}(x,y_j,Dw_j(x))\ -\ \theta\lambda_j\Big[w_j(x) + \Phi\big(-w_j(x),\,-w_{\bar\jmath}(x)\big)\Big],\quad j=1,2,\,\,\bar\jmath=3-j,
\end{aligned}
\end{equation}
subject to the state constraint $x\ge\ux$ and $r\ux+y_j-c_j(\ux)\ge0$.  The first-order condition for the maximization in \eqref{Eq: HJB-w} gives the optimal consumption
\begin{equation}\label{c-star}
c_j(x)=\left(\frac{Dw_j(x)}{\rho(\p^{-1}-1)}\right)^{-\p}.
\end{equation}
It is to be replaced by the constrained value $c_j(\ux)=r\ux+y_j$ whenever the unconstrained rule would violate the state constraint. The Hamiltonian in \cref{Eq: HJB-w} is given by 
\begin{equation}\label{Eq: H Dw}
  \mathcal{H}(x,y,p)=\left\{
      \begin{aligned}
     \quad (rx+y)p-\frac{1}{1-\psi}\left(\frac{\rho(1-\psi)}{\psi}\right)^\psi p^{1-\psi},\quad \quad & \text{if }\quad p\geq 0,\\
       +\infty, \quad \quad  & \text{if }\quad p< 0.
     \end{aligned}
   \right.
\end{equation}
The new Hamiltonian $\mathcal{H}(x,y,p)$ is convex w.r.t. $p$, and for fixed $(x,y)$ attains its minimum at 
$
p_{\min}(x,y)=\frac{\rho(1-\psi)}{\psi}(rx+y)^{-1/\psi}. 
$
\begin{remark}
When $\theta=1$ (i.e.\ $\gamma=\p^{-1}$), $\Phi(-w_j,-w_{\bar\jmath})=-w_{\bar\jmath}$, so the switching term reduces to $\lambda_j(w_{\bar\jmath}-w_j)$, recovering the HJB equations with CRRA utility in \cite{SL_HA}.
\end{remark}
We define sub- and supersolutions of the \cref{Eq: HJB-w} similarly as in \cref{def vis sol} and then prove the subsolution lies below the supersolution. We give only the first steps of the proof, as the only significant difference from the setting of \cite[Theorem 2.3]{SL_HA} is the presence of the nonlinear switching terms.

\begin{theorem}[Comparison principle]\label{Thm: comparison w}
Assume $0<\theta<1$, if $u=(u_1,u_2)$ and $w=(w_1,w_2)$ are sub- and supersolutions then
$u_j(x) \ \le\ w_j(x)$ for all $x\ge\ux,\ j=1,2$.
\end{theorem}

\begin{proof}
Suppose, for contradiction, that
$
\delta := \sup_{j\in\{1,2\},\,x\ge\ux}\big(u_j(x)-w_j(x)\big) \ >\ 0,
$
attained at some $(j^*,x^*)$ with $x^*\geq \ux$. By the maximality of $\delta$ (over both $j$ and $x$):
\begin{equation}\label{Eq: deltaj}
u_{j^*}(x^*)-w_{j^*}(x^*)=\delta, \qquad u_{\bar\jmath^*}(x^*)-w_{\bar\jmath^*}(x^*)\ \le\ \delta.
\end{equation}
From \cref{Eq: deltaj} and since $\Phi$ is increasing in its second argument, we infer
\beq\label{Eq: remainder-ineq}
\begin{aligned}
&\big[w_{j^*}(x^*)+\Phi(-w_{j^*}(x^*),-w_{\bar\jmath^*}(x^*))\big] -\big[u_{j^*}(x^*)+\Phi(-u_{j^*}(x^*),-u_{\bar\jmath^*}(x^*))\big] \\
= {}&-\delta+\Phi(-w_{j^*}(x^*),-w_{\bar\jmath^*}(x^*))-\Phi(-w_{j^*}(x^*)-\delta,-u_{\bar\jmath^*}(x^*))\\
\leq {}&-\delta+\Phi(-w_{j^*}(x^*),-w_{\bar\jmath^*}(x^*))-\Phi(-w_{j^*}(x^*)-\delta,-w_{\bar\jmath^*}(x^*)-\delta)
\leq  0,
\end{aligned}
\eeq
where the last inequality is yielded by \cref{lem:Phi-shift}. Multiplying \cref{Eq: remainder-ineq} by $\theta\lambda_{j^*}$  yields the relationship between the switching terms of the HJB equations:
$$
-\theta\lambda_{j^*}[u_{j^*}(x^*)+\Phi(-u_{j^*}(x^*),-u_{\bar\jmath^*}(x^*))]\ \le\ -\theta\lambda_{j^*}[w_{j^*}(x^*)+\Phi(-w_{j^*}(x^*),-w_{\bar\jmath^*}(x^*))].
$$
The rest of the proof proceeds similarly as \cite[Theorem 2.3]{SL_HA} and we omit the details. 
\end{proof}

\begin{theorem}
There exists a unique viscosity solution $v=(v_1,v_2)$ to the HJB system \cref{Eq: HJB}. 
\end{theorem}
\begin{proof}
We know from \cref{Prop: v sol early}, or by using Perron's method, that there exists $v=(v_1,v_2)$ which is a viscosity solution of \cref{Eq: HJB}. We infer there exists a viscosity solution $w$ of \cref{Eq: HJB-w} given by the change of variables \cref{Eq: v to w}, since this mapping is monotone increasing and differentiable. We know from \cref{Thm: comparison w} that $w$ is the unique solution of \cref{Eq: HJB-w}, hence $v$ is the unique solution of \cref{Eq: HJB}.
\end{proof}

\section{The semi-Lagrangian method}\label{Sec: SL 1}
\subsection{The SL scheme and first elements of analysis}\label{Subsec: discrete HJB}
Given a step size $\D x=(\bar{x}-\ux)/I$, $I\in\N$, let $x_i=\ux+i \Dx$, $i=0,\dots,I$, be the points of the space grid $\mathcal{G}^{\D x}=\{x_i\}_{i=0}^I$. We denote $\D=(h,\Dx)$ where $h$ is the time step. We denote by $\mathbb{I}[\cdot]$ the standard piecewise-linear
interpolation operator that maps a grid function
$\varphi=(\varphi_i)_{i=0}^I\in B(\mathcal{G}^{\D x})$ to a function on
$[\ux,\bar x]$. Specifically,
$
\mathbb{I}[\varphi](x)
=
\sum_{i=0}^I\varphi_i\beta_i(x),
$
where $(\beta_i)_{i=0}^I$ denotes the standard piecewise-linear
$\mathbb{Q}_1$ nodal basis associated with the grid. 
We recall that $\mathbb I$ is order preserving, i.e.,
\begin{equation}\label{Eq: I monotone}
\varphi\le \phi \ \text{ on }\ \mathcal G^{\D x}\quad\Longrightarrow\quad \mathbb I[\varphi](z)\le \mathbb I[\phi](z)\ \ \forall z\in[\ux,\bar x].
\end{equation}
We define the admissible set at a grid point $x_i$
\beq\label{Eq: C_j set}
\mC^{\D}_{j}(x_i)=[0,+\infty)\cap \left[\frac{x_i-\bar{x}}{h}+rx_i+y_j,\,\frac{x_i-\ux}{h}+rx_i+y_j\right],
\eeq
which enforces both the borrowing constraint at $\ux$ and the upper limit $\bar x$. 
Throughout the paper, we impose the standing assumption
$
0<h<\min\left\{
\frac{\theta}{\rho},
\frac{1}{\lambda_1},
\frac{1}{\lambda_2}
\right\}.
$
We now consider the semi-Lagrangian (SL) scheme obtained by a first order (in $h$) time discretization of the dynamic programming equation: the transport term is treated by evaluating the interpolated value $\mathbb I[V_j]$ at the point reached after one time step along the discrete characteristic, while the aggregator $f(c,v_j)$ is kept at the grid point itself. By using a Bernoulli process to approximate the Poisson process with rate $\lambda_j$, we construct the scheme as in \cite{SL_HA}, with $j\in \{1,2\}$ and $\bjm=3-j$: 
$$
\begin{aligned}
&V^{\D}_{i,j}=\sup_{c\in \mC^{\D}_{j}(x_i)}\left\{
hf\lc c,V^{\D}_{i,j}\rc +\left[\lambda_j hV^{\D}_{i,\bjm} 
+(1-\lambda_jh)\mathbb{I}[V^{\D}_j]\lc x_i+hs_{i,j}(c)\rc \right]
\right\},\\
& s_{i,j}(c)=rx_i+y_j-c,
\end{aligned}
$$
where $V^{\D}_{k,j}$ denotes the approximation of $v_j(x_k)$. We observe that with \cref{EZ flow} and  $(1+h\zeta)^{-1}=1-\zeta h+o(h)$, the above scheme is equivalent (up to first order in $h$) to 
\begin{equation}\label{Eq: HJBSL}
\begin{aligned}
V^{\D}_{i,j}=\sup_{c\in \mC^{\D}_{j}(x_i)}\left\{
h\mathcal{F}\lc c,V^{\D}_{i,j}\rc +(1-\zeta h)\left[\lambda_j hV^{\D}_{i,\bjm} 
+(1-\lambda_jh)\mathbb{I}[V^{\D}_j]\lc x_i+hs_{i,j}(c)\rc \right]
\right\}.
\end{aligned}
\end{equation}
We observe that if $\theta=1$ then $\rho=\zeta$ and $\mathcal{F}$ does not depend on $V$: we recover the SL scheme from \cite{SL_HA}. We use the notation $U\preceq V$ (resp. $U\succeq V$) introduced in \cref{Eq: I monotone}, $U_{i,j}\le V_{i,j}$ (resp. $U_{i,j}\ge V_{i,j}$) for all $i,j$.\par
We now introduce the discrete sub- and supersolutions of \cref{Eq: HJBSL}, analogous to \cref{def vis sol} at the continuous level.
\begin{definition}\label{Def: disc subsuper SL}
A grid function pair $U^{\D}=(U^{\D}_1,U^{\D}_2)$, $U^{\D}_j=(U^{\D}_{i,j})_{i=0}^I$, is a subsolution (resp. a supersolution) of \cref{Eq: HJBSL} if for every $i,j$, 
\begin{equation}\label{Eq: subsuper SL}
U^{\D}_{i,j}\ \le\ (\text{resp.}\ \ge)\ \sup_{c\in \mC^{\D}_{j}(x_i)}\left\{
h\mathcal{F}\lc c,U^{\D}_{i,j}\rc +\lc 1-\zeta h\rc \left[\lambda_j hU^{\D}_{i,\bjm} 
+(1-\lambda_jh)\mathbb{I}[U^{\D}_j]\lc x_i+hs_{i,j}(c)\rc \right]
\right\}.
\end{equation}
\end{definition}

We next give a discrete counterpart of \cref{Prop:sub-super}, providing barriers for the numerical solution. 
\begin{proposition}\label{Prop: disc barriers SL}
Let us define $\check\su,\check\sv$ as in \cref{Eq: sub super}. We set $\check U^{\D}_{i,j}=\check{\su}_j(x_i)$ and $\check V^{\D}_{i,j}=\check{\sv}_j(x_i)$. Then, the grid function pair $(\check U^{\D}_1,\check U^{\D}_2)$ is a subsolution of \cref{Eq: HJBSL}. Moreover, for sufficiently small $h>0$, the grid function pair $(\check V^{\D}_1,\check V^{\D}_2)$ is a supersolution of \cref{Eq: HJBSL}.
\end{proposition}
\begin{proof}
{\it Subsolution.} Let us fix a grid index $i$ and denote the null-saving consumption $\bar c_1(x_i)=rx_i+y_1$, which yields $s_{i,1}(\bar c_1(x_i))=0$. We can infer from \cref{Eq: sub super} and the scheme \cref{Eq: HJBSL} that at $j=1$,
$$
\begin{aligned}
&\sup_{c\in \mC^{\D}_{1}(x_i)}\left\{
h\mathcal{F}\lc c,\check U^{\D}_{i,1}\rc +(1-\zeta h)\left[\lambda_1h\check U^{\D}_{i,2} 
+(1-\lambda_1h)\mathbb{I}[\check U^{\D}_{1}]\lc x_i+hs_{i,1}(c)\rc \right]
\right\}-\check U^{\D}_{i,1}\\
\geq {}&h\mathcal{F}\lc \bar{c}_1(x_i),\check U^{\D}_{i,1}\rc-\frac{\rho h}{\theta}\check U^{\D}_{i,1}=0
\end{aligned}
$$
where for the last equality we used the algebraic identity $\mathcal F(\bar c_1(x),\check\su_1(x))=\frac{\rho}{\theta}\check\su_1(x)$. Therefore \cref{Eq: subsuper SL} holds at $j=1$. The same argument, with $\bar c_2(x_i)=rx_i+y_2$ and using $y_2>y_1$, $\mathcal F(\bar c_2(x_i),\check U^{\D}_{i,2})\ge \mathcal F(\bar c_1(x_i),\check U^{\D}_{i,2})=\frac{\rho}\theta \check U^{\D}_{i,2}$ (\cref{EZ flow} being increasing in $c$), yields \cref{Eq: subsuper SL} at $j=2$.\par 

{\it Supersolution.} Since $\check\sv_2$ is concave, we know that for every $c\in \mC^{\D}_{2}(x_i)$,
$$
\mathbb I[\check V^{\D}_2](x_i+hs_{i,2}(c))\ \le\ \check\sv_2(x_i+hs_{i,2}(c))\ \le\check V^{\D}_{i,2}+hs_{i,2}(c)D\check\sv_2(x_i).
$$
Since $\mathcal F\lc c,\check V^{\D}_{i,2}\rc<0$, $\mathcal F\lc c,\check V^{\D}_{i,2}\rc\leq (1-\lambda_2h)(1-\zeta h)\mathcal F\lc c,\check V^{\D}_{i,2}\rc$. 

\begin{align*}
&\sup_{c\in \mC^{\D}_{2}(x_i)}\left\{h\mathcal F\lc c,\check V^{\D}_{i,2}\rc+(1-\zeta h)\big[\lambda_2h\check V^{\D}_{i,1}+(1-\lambda_2h)\mathbb I[\check V^{\D}_2](x_i+hs_{i,2}(c))\big]\right\}-\check V^{\D}_{i,2}\\
\leq {}& h(1-\lambda_2h)(1-\zeta h)\lc \sup_{c\geq 0}\left\{\mathcal F\lc c,\check\sv_2(x_i)\rc+(rx_i+\bar{y}-c)D\check\sv_2(x_i)-\zeta \check\sv_2(x_i)\right\}\rc\\
&+h\left[\zeta \lc1-(1-\lambda_2h)(1-\zeta h)\rc (-\check\sv_2(x_i))+(1-\lambda_2h)(1-\zeta h)(y_2-\bar{y})D\check\sv_2(x_i)\right]\\
\leq {}&0.
\end{align*}
For the last inequality, we used $y_2-\bar{y}<0$ and $h$ sufficiently small. We have therefore obtained the $\geq$ inequality in \cref{Eq: subsuper SL} with $j=2$. The proof for $j=1$ is similar hence we omit it. 
\end{proof}

\subsection{A discrete comparison principle}\label{Sec: algo late}
 We now derive the comparison principle satisfied by the scheme \cref{Eq: HJBSL}. For the case $\theta=1$, a proof has been given in \cite[Proposition 4.6]{SL_HA}. For the case  $\theta>1$, we only need a straightforward extension of the same method.
We mainly focus on discussing the case $0<\theta<1$. As in \cref{Eq: v to w}, we perform the change of variable for the solution $V$ to \cref{Eq: HJBSL}:
\begin{equation}\label{Eq: V to W}
(1-\g)V_{j}= \lc-W_{j}\rc^{1/\theta}, \qquad j=1,2,
\end{equation}
A direct computation yields 
\begin{equation}\label{Eq: HJBSL-W}
\begin{aligned}
&W_{i,j}\\
={}&\ \sup_{c\in\mathcal C^\Delta_j(x_i)}\left\{-h\zeta\, c^{1-\p^{-1}} \ -\ (1-\zeta h)(1-\lambda_j h)\big(-W_{i,j}\big)^{1-1/\theta}\,\mathbb{I}\big[\big(-W_j\big)^{1/\theta}\big]\big(x_i+hs_{i,j}(c)\big)\right\}\\
{}& -(1-\zeta h)\,\lambda_j h\, \Phi\big(-W_{i,j},-W_{i,\bar\jmath}\big).
\end{aligned}
\end{equation}
The scheme \cref{Eq: HJBSL-W} in fact corresponds to the system of HJB equations \cref{Eq: HJB-w}. 
\begin{proposition}\label{Prop: comparison SL W}
Let us assume $0<\theta<1$. If $U$ and $W$ are respectively a sub- and a supersolution of \cref{Eq: HJBSL-W}, then $U\preceq W$.
\end{proposition}
\begin{proof}
Suppose, arguing by contradiction, that
\beq\label{Eq: U-W=delta}
U_{i^*,j^*}-W_{i^*,j^*}=\max_{i,j}\{U_{i,j}-W_{i,j}\}=\delta>0.
\eeq
Since $\mC^{\D}_{j^*}(x_{i^*})$ is compact and 
$$c\mapsto -h\zeta\, c^{1-\p^{-1}} -(1-\zeta h)(1-\lambda_{j^*} h)\big(-U_{i^*,j^*}\big)^{1-1/\theta}\,\mathbb{I}\big[\big(-U_{j^*}\big)^{1/\theta}\big]\big(x_{i^*}+hs_{i^*,j^*}(c)\big)$$
 is continuous, the supremum defining the subsolution property of $U$ at $(i^*,j^*)$ is attained at some $\widehat{c}\in \mC^{\D}_{j^*}(x_{i^*})$. Moreover, there exist two grids $x_k$ and $x_{k+1}$ such that $x_k\leq x_{i^*}+hs_{i^*,j^*}(\widehat{c})\leq x_{k+1}$. We denote by 
 $
 \beta:=\frac{x_{k+1}-x_{i^*}-hs_{i^*,j^*}(\widehat{c})}{\D x}.
 $
 We obtain
 $$
 \mathbb{I}\big[\big(-U_{j^*}\big)^{1/\theta}\big]\big(x_{i^*}+hs_{i^*,j^*}(\widehat{c})\big)=\beta \big(-U_{k,j^*}\big)^{1/\theta}+(1-\beta)\big(-U_{k+1,j^*}\big)^{1/\theta},
 $$
 which yields
\begin{equation}\label{Eq: HJBSL-W sub}
\begin{aligned}
U_{i^*,j^*}
\leq {}& -h\zeta\, \widehat{c}^{1-\p^{-1}} - (1-\zeta h)\lambda_{j^*} h\, \Phi\big(-U_{i^*,j^*},-U_{i^*,\bar\jmath^*}\big)\\
{}&-\ (1-\zeta h)(1-\lambda_{j^*}h)\lc \beta \Phi \big(-U_{i^*,j^*},-U_{k,j^*}\big)+(1-\beta)\Phi\big(-U_{i^*,j^*},-U_{k+1,j^*}\big)\rc.
\end{aligned}
\end{equation}
Since $W$ is a supersolution of \cref{Eq: HJBSL-W} and $\widehat{c}$ is an admissible control for 
$$
\sup_{c\in\mathcal C^\Delta_{j^*}(x_{i^*})}\left\{-h\zeta\, c^{1-\p^{-1}} \ -\ (1-\zeta h)(1-\lambda_{j^*} h)\big(-W_{i^*,j^*}\big)^{1-1/\theta}\,\mathbb{I}\big[\big(-W_{j^*}\big)^{1/\theta}\big]\big(x_{i^*}+hs_{i^*,j^*}(c)\big)\right\},
$$
we can infer
\begin{equation}\label{Eq: HJBSL-W super}
\begin{aligned}
W_{i^*,j^*}
\geq {}& -h\zeta\, \widehat{c}^{1-\p^{-1}} - (1-\zeta h)\lambda_{j^*} h\, \Phi\big(-W_{i^*,j^*},-W_{i^*,\bar\jmath^*}\big)\\
{}&-\ (1-\zeta h)(1-\lambda_{j^*}h)\lc \beta \Phi \big(-W_{i^*,j^*},-W_{k,j^*}\big)+(1-\beta)\Phi\big(-W_{i^*,j^*},-W_{k+1,j^*}\big)\rc.
\end{aligned}
\end{equation}
We observe that \cref{Eq: U-W=delta} yields
$$
U_{k,j^*}\leq W_{k,j^*}+\delta,\quad U_{k+1,j^*}\leq W_{k+1,j^*}+\delta,\quad U_{i^*,\bar\jmath^*}\leq W_{i^*,\bar\jmath^*}+\delta,
$$
hence we can deduce from \cref{lem:Phi-shift} that 
\beq\label{Eq: shift Phi W}
\begin{aligned}
&\Phi\big(-W_{i^*,j^*},-W_{i^*,\bar\jmath^*}\big)-\Phi\big(-U_{i^*,j^*},-U_{i^*,\bar\jmath^*}\big)\leq \delta,\\
&\Phi \big(-W_{i^*,j^*},-W_{k,j^*}\big)-\Phi \big(-U_{i^*,j^*},-U_{k,j^*}\big)\leq \delta,\\
& \Phi \big(-W_{i^*,j^*},-W_{k+1,j^*}\big)-\Phi \big(-U_{i^*,j^*},-U_{k+1,j^*}\big)\leq \delta.
\end{aligned}
\eeq
By comparing the two inequalities \cref{Eq: HJBSL-W sub} and \cref{Eq: HJBSL-W super}, while using \cref{Eq: shift Phi W}, we deduce $\delta\leq (1-\zeta h)\delta$ which is a contradiction. 
\end{proof}

\begin{proposition}\label{Prop: comparison SL}
If $U$ and $V$ are respectively a sub- and a supersolution of \cref{Eq: HJBSL} (\cref{Def: disc subsuper SL}), then $U\preceq V$.
\end{proposition}
\begin{proof}
The proof in the case $\theta=1$ is the same as \cite[Proposition 4.6]{SL_HA}. We now consider the other two cases.\par 
{\it{Case 1}}: $\theta>1$. Suppose, arguing by contradiction, that
$$
U_{i^*,j^*}-V_{i^*,j^*}=\max_{i,j}\{U_{i,j}-V_{i,j}\}=\delta>0.
$$
Since $\mC^{\D}_{j^*}(x_{i^*})$ is compact and $c\mapsto h\mathcal F(c,U_{i^*,j^*})+(1-\zeta h)(1-\lambda_{j^*}h)\mathbb I[U_{j^*}](x_{i^*}+hs_{i^*,j^*}(c))$ is continuous, the supremum defining the subsolution property of $U$ at $(i^*,j^*)$ is attained at some $\widehat{c}\in \mC^{\D}_{j^*}(x_{i^*})$, and, since $V$ is a supersolution, the same control $\widehat{c}$ is admissible in the corresponding inequality for $V$. Hence, with $\widehat{s}=s_{i^*,j^*}\lc \widehat{c}\rc$,
\begin{align*}
U_{i^*,j^*}&\le h\mathcal F(\widehat{c},U_{i^*,j^*})+(1-\zeta h)\big[\lambda_{j^*}hU_{i^*,\bjm^*}+(1-\lambda_{j^*}h)\mathbb I[U_{j^*}](x_{i^*}+h\widehat{s})\big],\\
V_{i^*,j^*}&\ge h\mathcal F(\widehat{c},V_{i^*,j^*})+(1-\zeta h)\big[\lambda_{j^*}hV_{i^*,\bjm^*}+(1-\lambda_{j^*}h)\mathbb I[V_{j^*}](x_{i^*}+h\widehat{s})\big].
\end{align*}
Since $U_{i^*,j^*}>V_{i^*,j^*}$, we know from $\theta>1$ and \cref{F monotone} that 
$$
\mathcal F(\widehat{c},U_{i^*,j^*})-\mathcal F(\widehat{c},V_{i^*,j^*})< 0.
$$
Since $U_{i^*,\bjm^*}-V_{i^*,\bjm^*}\le \delta$ and, by \cref{Eq: I monotone} together with $U_{j^*}-V_{j^*}\le \delta$ on the grid,
$$
\mathbb I[U_{j^*}](x_{i^*}+h\widehat{s})-\mathbb I[V_{j^*}](x_{i^*}+h\widehat{s})\le \delta.
$$
Subtracting the two inequalities above and combining with these two bounds, we obtain
$$
\delta=U_{i^*,j^*}-V_{i^*,j^*}\le (1-\zeta h)\big[\lambda_{j^*}h\,\delta+(1-\lambda_{j^*}h)\delta\big]=(1-\zeta h)\delta,
$$
i.e. $\zeta h\delta\le 0$, yielding a contradiction.\par 
{\it{Case 2}}: $0<\theta<1$. Let us denote 
$$
(1-\g)U= \lc-\widetilde{U}\rc^{1/\theta},\quad (1-\g)V= \lc-\widetilde{V}\rc^{1/\theta},
$$
then $\widetilde{U}$ and $\widetilde{V}$ are respectively sub- and supersolution of \cref{Eq: HJBSL-W}. From \cref{Prop: comparison SL W}, we know $\widetilde{U}\preceq \widetilde{V}$, hence $U\preceq V$. 
\end{proof}
By \cref{Prop: disc barriers SL}, $(\check U_1,\check U_2)$ and $(\check V_1,\check V_2)$ are respectively a subsolution and a supersolution of \cref{Eq: HJBSL}, so that \cref{Prop: comparison SL} yields the following discrete barrier property.
\begin{proposition}\label{Cor: sol bound SL}
The solution $V$ to \cref{Eq: HJBSL}, if it exists, is unique and satisfies $\check U\preceq V\preceq \check V$.
\end{proposition}
 \subsection{A Howard--Newton algorithm}
Following a similar design as in the finite-difference case (cf. \cite[Algorithm 4.1]{achdou2026finite}), \cref{Eq: HJBSL} can be solved by the following Howard--Newton algorithm, in which the policy evaluation step requires an inner loop of Newton iterations because the aggregator $\mathcal F(c,V_{i,j})$ depends nonlinearly on $V_{i,j}$ if $\theta \neq 1$.

\begin{algorithm}[h!]
\caption{Howard--Newton algorithm}
\label{Algo: Howard Newton}
\begin{algorithmic}[1]
\State Initialize $c^{(0)}_{i,j}=rx_i+y_j$, and set $n=0$
\Repeat \Comment{Outer loop}
\State \textbf{Saving policy:} $s^{(n)}_{i,j}=rx_i+y_j-c^{(n)}_{i,j}$
\State \textbf{Policy evaluation:} solve for $V^{(n)}_{i,j}$, $j=1,2$, $i=0,\dots,I$,
\beq\label{Eq: Howard solve V SL}
V^{(n)}_{i,j}=h\mathcal F\lc c^{(n)}_{i,j},V^{(n)}_{i,j}\rc +(1-\zeta h)\left[\lambda_j hV^{(n)}_{i,\bjm} 
+(1-\lambda_jh)\mathbb{I}[V^{(n)}_j]\lc x_i+hs^{(n)}_{i,j}\rc \right],
\eeq
by Newton's method (the map $V_{i,j}\mapsto h\mathcal F(c^{(n)}_{i,j},V_{i,j})$ being smooth in $V_{i,j}$)\\
\Comment{Inner loop}
\State \textbf{Consumption policy update:} 
\beq\label{Eq: Howard-Newton P-update SL}
c^{(n+1)}_{i,j}=\mathop{\arg\max}\limits_{c\in \mC^{\D}_{j}(x_i)} \left\{h\mathcal F\lc c,V^{(n)}_{i,j}\rc +(1-\zeta h)(1-\lambda_jh)\mathbb{I}[V^{(n)}_j]\lc x_i+hs_{i,j}(c)\rc\right\}
\eeq
\State $n\gets n+1$
\Until{$\sum_j\max_i\left\vert c^{(n)}_{i,j}-c^{(n-1)}_{i,j}\right\vert<\text{Tol}_c$}
\end{algorithmic}
\end{algorithm}

\begin{proposition}\label{Prop: convergence Howard}
We make the same assumption as \cref{Prop: disc barriers SL}. The sequence $V^{(n)}$ generated by \cref{Algo: Howard Newton} converges to a bounded grid function $V$, which is the unique solution to \cref{Eq: HJBSL}. 
\end{proposition}
\begin{proof}
From \cref{Eq: Howard solve V SL}, we observe that $V^{(n+1)}$ is the solution to the discrete system
\beq\label{Eq: SL linear n}
V_{i,j}=h\mathcal F\lc c^{(n+1)}_{i,j},V_{i,j}\rc +(1-\zeta h)\left[\lambda_j hV_{i,\bjm} 
+(1-\lambda_jh)\mathbb{I}[V_j]\lc x_i+hs^{(n+1)}_{i,j}\rc \right].
\eeq
From \cref{Eq: Howard solve V SL} and \cref{Eq: Howard-Newton P-update SL} at step $n$, 
$$
V^{(n)}_{i,j}\leq h\mathcal F\lc c^{(n+1)}_{i,j},V^{(n)}_{i,j}\rc +(1-\zeta h)\left[\lambda_j hV^{(n)}_{i,\bjm} 
+(1-\lambda_jh)\mathbb{I}[V^{(n)}_j]\lc x_i+hs^{(n+1)}_{i,j}\rc \right],
$$
i.e. $V^{(n)}$ is a subsolution to \cref{Eq: SL linear n}. We can then use a similar discrete comparison principle as \cref{Prop: comparison SL} to show $V^{(n)}\preceq V^{(n+1)}$. Since $V^{(n)}$ is a subsolution to \cref{Eq: HJBSL}, we infer from \cref{Prop: comparison SL} that for all $n$, $V^{(n)}\preceq \check{V}^\D$. Therefore $V^{(n)}$ converges to a limit $V$, and it is the unique solution of \cref{Eq: HJBSL} from \cref{Cor: sol bound SL}. 
\end{proof}
We now discuss the implementation of the consumption policy update \cref{Eq: Howard-Newton P-update SL}.
For a generic pair of grid functions $V=(V_1,V_2)$ we define the finite-difference operator
\beq\label{Eq: discrete D}
\begin{aligned}
&{\mathbf{D}}^+(V_{i,j})=\frac{V_{i+1,j}-V_{i,j}}{\D x},\,\,0\leq i\leq I-1\quad {\text{and}}\quad {\mathbf{D}}^-(V_{i,j})=\frac{V_{i,j}-V_{i-1,j}}{\D x},\,\,1\leq i\leq I,\\
&{\mathbf{D}}(V_{i,j})=
\left\{\begin{aligned}
\quad &{\mathbf{D}}^+(V_{0,j}),\quad &&i=0,\\
&\frac{V_{i+1,j}-V_{i-1,j}}{2\D x},&&0<i<I,\\
&{\mathbf{D}}^-(V_{I,j}),&&i=I.
\end{aligned}\right.
\end{aligned}
\eeq
We observe that $\mathbf{D}^+(V_{k,j})$ is the slope of the interpolant $\mathbb I[V_j]$ on the cell $[x_k,x_{k+1}]$. For a grid point $(x_i,y_j)$ and an arrival point $z$, let
\beq\label{Eq: chi}
\chi_{i,j}(z)=rx_i+y_j+\frac{x_i-z}{h}
\eeq
be the consumption at $(x_i,y_j)$ that leads to $z=x_i+hs_{i,j}(c)$. We consider the problem to solve for the policy update step \cref{Eq: Howard-Newton P-update SL} in \cref{Algo: Howard Newton}:
\beq\label{Eq: c hat}
\widehat{c}_{i,j}:=\mathop{\arg\max}\limits_{c\in \mC_{j}^\D(x_i)} \left\{h\mathcal F\lc c,V_{i,j}\rc +(1-\zeta h)(1-\lambda_jh)\,\mathbb{I}[V_j]\lc x_i+hs_{i,j}(c)\rc\right\},
\eeq
with $V=V^{(n)}$ in the policy update. For each cell $0\leq k\leq I-1$ let
\beq\label{Eq: set Ck}
\mC^{k}_{j}(x_i):=\left\{c>0:\ x_i+hs_{i,j}(c)\in [x_k,x_{k+1}]\right\}
=\left\{c>0:\ \chi_{i,j}(x_{k+1})\leq c\leq \chi_{i,j}(x_k)\right\},
\eeq
and
\beq\label{Eq: c hat k argmax}
\widehat{c}_{i,j}^k:=\mathop{\arg\max}\limits_{c\in \mC^{k}_{j}(x_i)} \left\{\underbrace{h\mathcal F\lc c,V_{i,j}\rc +(1-\zeta h)(1-\lambda_jh)\,\mathbb{I}[V_j]\lc x_i+hs_{i,j}(c)\rc}_{\Phi_{i,j}(c)} \right\}.
\eeq
Since the cells cover $[x_0,x_I]$, we have $\mC^{\D}_j(x_i)=\bigcup_{k}\mC^k_j(x_i)$.  
\begin{lemma}\label{Lem: c hat k}
Fix $(i,j)$ and a cell $k$ with $\chi_{i,j}(x_k)>0$ (otherwise $\mC^k_j(x_i)=\emptyset$). If $\mathbf{D}^+V_{k,j}>0$, let
\beq\label{Eq: c k}
c_{i,j}^k=\lc\frac{\rho}{(1-\zeta h)(1-\lambda_jh)\,\mathbf{D}^+V_{k,j}}\rc^{\p}\lc(1-\g)V_{i,j}\rc^{\frac{1-\g \p}{1-\g}}.
\eeq
Then the maximizer $\widehat{c}_{i,j}^k$ in \cref{Eq: c hat k argmax} exists, is unique, and is given by
\beq\label{Eq: c hat k}
\widehat{c}_{i,j}^k=\left\{
      \begin{aligned}
     \quad &\max\Bigg\{\min\left\{c_{i,j}^k,\ \chi_{i,j}(x_k)\right\},\ \chi_{i,j}(x_{k+1})\Bigg\}\ &&{\text{if}}\ \mathbf{D}^+V_{k,j}>0,\\
     \quad &\chi_{i,j}(x_k) &&{\text{if}}\ \mathbf{D}^+V_{k,j}\leq 0.
      \end{aligned}
   \right.
\eeq
\end{lemma}
\begin{proof}
For $c\in\mC^k_j(x_i)$ the arrival point $x_i+hs_{i,j}(c)=x_i+h(rx_i+y_j-c)$ lies in $[x_k,x_{k+1}]$, where
$\mathbb I[V_j](z)=V_{k,j}+\mathbf D^+V_{k,j}(z-x_k)$. The objective $\Phi_{i,j}(c)$ in \cref{Eq: c hat k argmax} is strictly concave in $c$, with derivative $h\,[\mathcal F_c(c,V_{i,j})-(1-\zeta h)(1-\lambda_jh)\mathbf D^+V_{k,j}]$.
If $\mathbf D^+V_{k,j}>0$, the
maximizer on the interval $[\chi_{i,j}(x_{k+1}),\chi_{i,j}(x_k)]\cap(0,\infty)$ is $c^k_{i,j}$ projected onto it. Since $c^k_{i,j}>0$ and $\chi_{i,j}(x_k)>0$, the projection is
given by the first line of \cref{Eq: c hat k}. If $\mathbf D^+V_{k,j}\leq0$, then $\mathcal F_c(c,V_{i,j})-(1-\zeta h)(1-\lambda_jh)\mathbf D^+V_{k,j}>0$, hence the objective $\Phi_{i,j}(c)$ is increasing in $c$ and the maximizer is the largest admissible consumption $\chi_{i,j}(x_k)$.
\end{proof}
 
By \cref{Lem: c hat k} and $\mC_j(x_i)=\bigcup_k\mC^k_j(x_i)$, the global maximizer in \cref{Eq: c hat} is obtained by comparing the cell-wise maximizers:
\beq\label{Eq: c hat argmax k}
\widehat{c}_{i,j}=\widehat{c}_{i,j}^{\,k^\ast},\qquad
k^\ast\in\mathop{\arg\max}\limits_{k\in K_{i,j}}\left\{h\mathcal F\lc \widehat{c}_{i,j}^k,V_{i,j}\rc +(1-\zeta h)(1-\lambda_jh)\,\mathbb{I}[V_j]\lc x_i+hs_{i,j}(\widehat{c}_{i,j}^k)\rc \right\},
\eeq
where $K_{i,j}=\{k:\chi_{i,j}(x_k)>0\}$ is the set of cells with admissible consumption, and ties are broken in favor of the smallest $k$. Computing \cref{Eq: c hat argmax k} for all $(i,j)$ solves \cref{Eq: c hat} exactly, but requires $O(I)$ operations per grid point and hence $O(I^2)$ per policy update. In \cref{App: argmax} we give an accelerated exact method that exploits the concavity of $V_j$ on groups of cells, and in \cref{App: argmax FOC} an approximate alternative based on the first-order condition.

\subsection{A monotone value iteration algorithm in the case \texorpdfstring{$0<\theta<1$}{0 < theta < 1}}\label{Sec: algo early}
Throughout this section we assume $0<\theta<1$. To solve the scheme \cref{Eq: HJBSL} in this case, we employ a direct value iteration algorithm, in which the value function is updated successively through an order preserving mapping. \par
We define the map $V=(V_1,V_2)=\Gamma(\widetilde{V}_1,\widetilde{V}_2)$: 
\begin{equation}\label{Eq: PsiSL}
\begin{aligned}
&V_{i,j}=\sup_{c\in \mC^{\D}_{j}(x_i)}\left\{
h\mathcal{F}\lc c,\widetilde{V}_{i,j}\rc +(1-\zeta h)\left[\lambda_j h\widetilde{V}_{i,\bjm} 
+(1-\lambda_jh)\mathbb{I}[\widetilde{V}_j]\lc x_i+hs_{i,j}(c)\rc \right]
\right\},
\end{aligned}
\end{equation}
 for $j=1,2$ and $i=0,\dots,I$.\par
\cref{Prop: disc barriers SL} yields directly the following result. 
\begin{lemma}\label{Lemma: early barrier SL}
For sufficiently small $h>0$ as in \cref{Prop: disc barriers SL}, we have the barrier properties
\begin{itemize}
\item[(i).] If $U=\Gamma(\check U_1,\check U_2)$, then $U\succeq \check U$.
\item[(ii).] If $V=\Gamma(\check V_1,\check V_2)$, then $V\preceq \check V$.
\end{itemize}
\end{lemma}

\begin{proposition}\label{Prop: Gamma}
Assume $0<\theta<1$. We have the following properties of $\Gamma$.
\begin{itemize}
\item Monotonicity. If $U=\Gamma(\widetilde U_1,\widetilde U_2)$, $V=\Gamma(\widetilde V_1,\widetilde V_2)$ and $\widetilde U\preceq \widetilde V$, then $U\preceq V$.
\item Invariance. If $\check U\preceq \widetilde V\preceq \check V$ and $V=\Gamma(\widetilde V_1,\widetilde V_2)$, then $\check U\preceq V\preceq \check V$.
\end{itemize}
\end{proposition}
\begin{proof}
From $0<\theta<1$ and \cref{F monotone}, we know that $\widetilde U\preceq \widetilde V$ yields $\mathcal{F}\lc c,\widetilde{U}_{i,j}\rc\leq \mathcal{F}\lc c,\widetilde{V}_{i,j}\rc$ for all $(i,j)$ and $c\in \mC^{\D}_{j}(x_i)$. Therefore, 
$$
\begin{aligned}
&\sup_{c\in \mC^{\D}_{j}(x_i)}\left\{
h\mathcal{F}\lc c,\widetilde{U}_{i,j}\rc +(1-\zeta h)\left[\lambda_j h\widetilde{U}_{i,\bjm} 
+(1-\lambda_jh)\mathbb{I}[\widetilde{U}_j]\lc x_i+hs_{i,j}(c)\rc \right]
\right\}\\
\leq {}&\sup_{c\in \mC^{\D}_{j}(x_i)}\left\{
h\mathcal{F}\lc c,\widetilde{V}_{i,j}\rc +(1-\zeta h)\left[\lambda_j h\widetilde{V}_{i,\bjm} 
+(1-\lambda_jh)\mathbb{I}[\widetilde{V}_j]\lc x_i+hs_{i,j}(c)\rc \right]
\right\}.
\end{aligned}
$$
Invariance follows from monotonicity and \cref{Lemma: early barrier SL}.
\end{proof}
Taking advantage of the monotonicity and invariance properties of $\Gamma$ established in \cref{Prop: Gamma}, we propose \cref{Algo: up}, initialized at the discrete counterpart of the subsolution $\check\su$, namely $\check U_{i,j}=\check{\su}_j(x_i)$.

\begin{algorithm}[h!]
\caption{Monotone value iteration algorithm}
\label{Algo: up}
\begin{algorithmic}[1]
\State Initialize $V^{(0)}_{i,j}=\check U_{i,j}$, and set $k=0$ 
\Repeat
\State $V^{(k+1)}=\Gamma \lc V^{(k)}_1,V^{(k)}_2\rc$
\State $k\gets k+1$
\Until{$\sum_j\max_i\dfrac{\left\vert V^{(k)}_{i,j}-V^{(k-1)}_{i,j}\right\vert}{1+\left\vert V^{(k-1)}_{i,j}\right\vert}<\text{Tol}_V$}
\end{algorithmic}
\end{algorithm}
\begin{theorem}\label{Thm: SL_VI convergence}
The sequence $V^{(k)}$ generated by \cref{Algo: up} converges to a bounded grid function $V$, which is the unique solution to \cref{Eq: HJBSL}. 
\end{theorem}
\begin{proof}
We know from \cref{Lemma: early barrier SL} that $V^{(0)} \preceq V^{(1)}$. From \cref{Prop: Gamma} and induction we deduce $V^{(k)} \preceq V^{(k+1)}$. Moreover, \cref{Prop: Gamma} shows that $V^{(k)} \preceq \check{V}$. Therefore, $V^{(k)}$ converges to a grid function $V$ such that $V=\Gamma(V)$. From \cref{Cor: sol bound SL}, $V$ is the unique solution to \cref{Eq: HJBSL}. 
\end{proof}
\begin{remark}
\cref{Algo: up} involves a maximization problem. In practice, we use \cref{Alg: Howard control update} to solve \cref{Eq: c hat} with $V_{i,j}=V_{i,j}^{(k)}$. This yields an optimal consumption $c^{(k+1)}$ and $\Gamma \lc V^{(k)}_1,V^{(k)}_2\rc$ is given by 
$
V^{(k+1)}_{i,j}=h\mathcal F\lc c^{(k+1)}_{i,j},V^{(k)}_{i,j}\rc +(1-\zeta h)\left[\lambda_j hV^{(k)}_{i,\bjm} 
+(1-\lambda_jh)\mathbb{I}[V^{(k)}_j]\lc x_i+hs^{(k+1)}_{i,j}\rc \right].
$
\end{remark}
\subsection{First order rate of convergence}
Next, we analyze the convergence rate of the solution to the scheme (on the grid) towards the solution of the HJB equation. Such result for general HJB equations has been considered in e.g. \cite{carlini2020semi}. The key ingredients in our analysis are the local regularity $v_j\in W^{2,\infty}_{\mathrm{loc}}(\ux,+\infty)$ (cf. \cref{Prop: v sol}) and a delicate expansion of $s_1$ and $Dv_1$ near the boundary $\ux$ established in \cref{Prop: Dv1 ux}. These estimates allow us to control the local interpolation error near $\ux$ and establish a first-order convergence rate. An analogous result for the finite difference scheme was established in \cite[Theorem 4.10]{achdou2026finite}. We start with preliminary results. 
\begin{lemma}\label{Lemma: interp}
We make the same assumptions as \cref{Prop: saving}, and $\bar{x}\in \mathcal{G}^{\D x}$ is sufficiently large such that $s_2(\bar{x})<0$. Let $v_j$ be the solution to the HJB equation \cref{HJB} and $s_j$ denotes the optimal saving policy. For any $x\in \mathcal{G}^{\D x}$, if $h=O(\D x)$ then there exists a constant $C>0$ such that
\beq\label{Eq: interp1}
-C\D x \leq \frac{\mathbb{I}[v_{j}]\lc x+hs_j(x)\rc-v_{j}(x+hs_j(x))}{h}\leq 0,
\eeq
and
\beq\label{Eq: interp2}
\frac{\mathbb{I}[v_{j}]\lc x+hs_j(x)\rc-v_{j}(x)}{h}-s_j(x)Dv_j(x)= O(\D x).
\eeq
\end{lemma}
\begin{proof}
From \cref{Prop: v sol}, $v_j$ is concave, hence $\mathbb{I}[v_{j}]\lc x+hs_j(x)\rc-v_{j}(x+hs_j(x))\leq 0$. Now we consider the lower bound in \cref{Eq: interp1} in two steps.\par 
{\it{Step 1.}} We consider the case $j=2$, which is the easier one. From \cref{Prop: saving}, we know $s_2(\ux)>0$. Hence by differentiating \cref{HJB} with $j=2$, we infer $D^2v_2(\ux)$ is bounded. Therefore, with $v_2\in W_{loc}^{2,\infty}$ we obtain $v_2\in W^{2,\infty}$. Let us denote $z_2=x+hs_2(x)$ for $x\in \mathcal{G}^{\D x}$, then there exists $x_k,\,x_{k+1}\in \mathcal{G}^{\D x}$ such that $x_k\leq z_2\leq x_{k+1}$. We obtain 
$$
\begin{aligned}
0\leq (x_{k+1}-z_2)(z_2-x_k)\leq \Delta x\min\{z_2-x_k,x_{k+1}-z_2\}\leq h|s_2(x)|\Delta x,
\end{aligned}
$$
hence
$$\begin{aligned}
\vert \mathbb{I}[v_{2}](z_2)-v_2(z_2)\vert\leq {}&\frac{1}{2}\|D^2v_2\|_{L^\infty}(x_{k+1}-z_2)(z_2-x_k)\\
\leq {}&\frac{1}{2}\|D^2v_2\|_{L^\infty}\|s_2\|_{L^\infty}h\D x,
\end{aligned}
$$
and this yields $$\sup_{x\in \mathcal{G}^{\D x}}\vert \mathbb{I}[v_{2}]\lc x+hs_2(x)\rc-v_{2}(x+hs_2(x))\vert \leq Ch\D x.$$\par
{\it{Step 2.}} We consider the case $j=1$, which is the more difficult case since $D^2v_1(x)$ blows up as $x\to \ux$. From \cref{Prop: Dv1 ux}, there exists a sufficiently small constant $\eta>0$ such that 
\beq\label{Eq: v1 eta bound}
 D^2v_1(x)\geq -\frac{C}{\sqrt{x-\ux}},\quad 0\leq -s_1(x)\leq C\sqrt{x-\ux},\quad -hs_1(x)\leq \D x\quad \forall x\in (\ux,\eta),
\eeq
where $C$ is a positive constant which does not depend on $x$, $h$ or $\D x$. The CFL condition $-hs_1(x)\leq \D x$ in \cref{Eq: v1 eta bound} is deduced from $-s_1(x)\leq C\sqrt{x-\ux}$ and $h=O(\D x)$. \par
Since $s_1(x)\leq 0$ for all $x\in [\ux,\bar{x}]$, we know $x+hs_1(x)$ is bracketed by two grids $x_k,\,x_{k+1}\in \mathcal{G}^{\D x}$ such that
\beq\label{Eq: point z}
z=x+hs_1(x)<x,\quad x_k\leq  z< x_{k+1}\leq x.
\eeq
 It follows from $v_1\in W_{loc}^{2,\infty}$ that $v_1\in W^{2,\infty}$ on the compact domain $[\eta',\bar{x}]$ where $\eta'=\ux+(\eta-\ux)/2$. From $h=O(\D x)$ and the boundedness of $s_1$ on $[\ux,\bar{x}]$, we know that $x_k\geq z-\D x>\eta'$ if $\D x$ is sufficiently small. Similarly as in {\it{Step 1}}, we can deduce 
$$
\sup_{x\in \mathcal{G}^{\D x}\cap [\eta,\bar{x}]}\vert \mathbb{I}[v_{1}]\lc x+hs_1(x)\rc-v_{1}(x+hs_1(x))\vert \leq \frac{1}{2}\|D^2v_1\|_{L^\infty[\eta',\bar{x}]}\|s_1\|_{L^\infty}h\D x.
$$\par
We now focus on the case $\ux<x\leq \eta$. From $s_1(x)<0$, \cref{Eq: v1 eta bound} and \cref{Eq: point z} we infer that in this regime
$$x-\D x= x_k\leq z<x_{k+1}=x\leq \eta.$$ 
The interpolation error can be decomposed as
\beq\label{v1 interp error}
\mathbb{I}[v_{1}](z)-v_1(z)
=\underbrace{\frac{x-z}{\D x}\int_{x_k}^z(\tau-x_k)D^2v_1(\tau)d\tau}_{(I)}+\underbrace{\frac{z-x_k}{\D x}\int_z^{x}(x-\tau)D^2v_1(\tau)d\tau}_{(II)}.
\eeq \par
Since $\tau-x_k\leq z-x_k\leq \D x$, $x-z=-hs_1(x)$ and $D^2v_1(\tau)<0$ in $(I)$, we infer
$$
(I)\geq -hs_1(x)\int_{x_k}^zD^2v_1(\tau)d\tau.
$$
Similarly $x-\tau\leq x-z=-hs_1(x)$ in $(II)$, and we deduce 
$$
(II)\geq -hs_1(x)\int_z^{x}D^2v_1(\tau)d\tau.
$$
Therefore, we infer from \cref{Eq: v1 eta bound} and \cref{v1 interp error} that there exists $C>0$ such that
\beq
\begin{aligned}
\mathbb{I}[v_{1}](z)-v_1(z)\geq{}& -hs_1(x)\int_{x_k}^{x}D^2v_1(\tau)d\tau\\
\geq {}&Ch\sqrt{x-\ux}\int_{x_k}^{x}\frac{-C}{\sqrt{\tau-\ux}}d\tau\\
\geq {}&-2C^2h\sqrt{x-\ux}(\sqrt{x-\ux}-\sqrt{x_k-\ux})\\
\geq {}&-2C^2h\sqrt{x-\ux}\frac{\D x}{\sqrt{x-\ux}+\sqrt{x_k-\ux}}\geq -2C^2h\D x. 
\end{aligned}
\eeq
Finally, we observe that $\mathbb{I}[v_{1}](\ux)-v_1(\ux)=0$ since $s_1(\ux)=0$. We have therefore completed the proof of \cref{Eq: interp1} with {\it{Step 1}} and {\it{Step 2}}. From \cref{Prop: v sol}, \cref{Prop: Dv1 ux} and the uniform boundedness of $s_j$, we infer
$$
\frac{v_{j}\lc x+hs_j(x)\rc-v_{j}(x)}{h}=s_j(x)Dv_j(x)+O(h),\quad \forall x\in [\ux,\bar{x}].
$$
We can then deduce \cref{Eq: interp2} from \cref{Eq: interp1}, $h=O(\D x)$ and 
$$
\frac{\mathbb{I}[v_{j}]\lc x+hs_j(x)\rc-v_{j}(x)}{h}=\frac{\mathbb{I}[v_{j}]\lc x+hs_j(x)\rc-v_{j}\lc x+hs_j(x)\rc}{h}+\frac{v_{j}\lc x+hs_j(x)\rc-v_{j}(x)}{h}.
$$
\end{proof}

\begin{lemma}\label{Lemma: h=O(Dx) c}
We make the same assumptions as in \cref{Lemma: interp}. Let $(v_1,v_2)$ be the solution to \cref{HJB} and $\widehat{c}_j(x)$ be the optimal control such that, for all $x_i\in \mathcal{G}^{\D x}$,
\beq\label{Eq: def widehat c}
H(x_i,y_j,v_j(x_i),Dv_j(x_i))=\mathcal{F}(\widehat{c}_j(x_i),v_j(x_i)) +(\underbrace{rx_i +y_j- \widehat{c}_j(x_i)}_{s_j(x_i)})Dv_j(x_i),\quad j\in \{1,2\}.
\eeq
 Suppose $h=O(\D x)$, then with sufficiently small $\D x$ we have $\widehat{c}_j(x_i)\in  \mC^{\D}_{j}(x_i)$. 
\end{lemma}
\begin{proof}
We first prove the case $j=1$. First we observe from \cref{Prop: saving} that $s_1(\ux)=0$, hence $\widehat{c}_1(\ux)\in  \mC^{\D}_{1}(\ux)$. 
We now consider the case $x_i>\ux$, hence $x_i-\ux \geq \D x$ and $\D x\leq \sqrt{\D x}\sqrt{x_i-\ux}$. From \cref{Eq: -s1 bound}, there exists a constant $\eta>\ux$ such that $\vert s_1(x_i)\vert =O(\sqrt{x_i-\ux})$ for all $x_i\in (\ux,\eta]$.
Therefore for all $x_i\in (\ux,\eta]$, 
$
h\vert s_1(x_i)\vert =O\lc h\sqrt{x_i-\ux}\rc=O\lc \D x\sqrt{x_i-\ux}\rc=O\lc \sqrt{\D x}(x_i-\ux)\rc.
$
Moreover, with sufficiently small $\D x$ and $h=O(\D x)$, we can infer from the boundedness of $s_1$ that 
$$
\sup_{x_i\in \mathcal{G}^{\D x},\,\eta<x_i\leq \bar{x}}h\vert s_1(x_i)\vert \leq \eta-\ux\leq x_i-\ux.
$$
This implies $x_i+hs_1(x_i)\geq\ux$ holds for all $x_i\in \mathcal{G}^{\D x}$ if $\D x$ is sufficiently small. The fact $x_i+hs_1(x_i)\leq \bar{x}$ is trivial since $s_1(x_i)<0$ for all $x>\ux$. We conclude $\widehat{c}_1(x_i)\in  \mC^{\D}_{1}(x_i)$.\par
The case with $j=2$ is simpler, since there exists $z\leq \widehat{x}$ such that $s_2\geq 0$ on $[\ux,z]$ and $s_2\leq 0$ on $[\widehat{x},\bar{x}]$. Therefore we omit the details. 
\end{proof}

\begin{theorem}\label{Thm: convergence rate}
We make the same assumptions as in \cref{Prop: saving}, and $\bar{x}\in \mathcal{G}^{\D x}$ is sufficiently large such that $s_2(\bar{x})<0$. Then for sufficiently small $\D x$ and $h=O(\D x)$, there exists a constant $C>0$, independent of $\D x$, such that
\beq\label{Eq: error}
\max_{j,\,x\in \mathcal{G}^{\D x}} \left\vert v_j(x)-V^{\D}_{j}(x)\right\vert \leq C\D x.
\eeq
\end{theorem}
\begin{proof} In this proof, $C$ denotes a generic positive constant that may change from line to line. \par
{\it{Step 1.}}
Let us define 
\beq\label{Eq: def sigma}
\sigma=\max_{j,\, x\in \mathcal{G}^{\D x}} \{v_j(x)-V^{\D}_{j}(x)\}=v_{j^*}(x^*)-V^{\D}_{j^*}(x^*),
\eeq
and suppose $\sigma>0$. Since $v_j\in C^1(\ux,+\infty)$, we know that if $x^*>\ux$ then
\beq\label{Eq: v(x) sub}
\frac{\rho}{\theta} v_{j^*}(x^*)=
	H\lc x^*,y_{j^*},v_{j^*}(x^*),Dv_{j^*}(x^*)\rc+\lambda_{j^*}(v_{\bar \jmath^*}(x^*)-v_{j^*}(x^*)).
\eeq
Since $Dv_j$ is uniformly continuous in $[\ux,+\infty)$, \cref{Eq: v(x) sub} holds if $x^*=\ux$. \par 
We denote by $\widehat{c}_{j^*}(x^*)$ the optimal control given by \cref{Eq: def widehat c}, and $\widehat{s}_{j^*}(x^*)=rx^* +y_{j^*}- \widehat{c}_{j^*}(x^*)$. We know from \cref{Eq: def sigma} that for all $j$ and $x\in \mathcal{G}^{\D x}$, $V^{\D}_j(x)\geq v_j(x)-\sigma$, hence
$$
V^{\D}_{\bar \jmath^*}(x^*)\geq v_{\bar \jmath^*}(x^*)-\sigma,\quad \mathbb{I}[V^{\D}_{j^*}]\lc x^*+h\widehat{s}_{j^*}(x^*)\rc\geq \mathbb{I}[v_{j^*}]\lc x^*+h\widehat{s}_{j^*}(x^*)\rc-\sigma.
$$
From $v_{j^*}(x^*)-\sigma=V^{\D}_{j^*}(x^*)$ and the SL scheme \cref{Eq: HJBSL}, we infer
$$
\begin{aligned}
&v_{j^*}(x^*)-\sigma \\
\geq{}& \sup_{c\in \mC^{\D}_{j^*}(x^*)}\left\{
h\mathcal{F}\lc c,v_{j^*}(x^*)\rc +(1-\zeta h)\left(\lambda_{j^*}hv_{\bar \jmath^*}(x^*)
+(1-\lambda_{j^*}h)\mathbb{I}[v_{j^*}]\lc x^*+h(rx^*+y_{j^*}-c)\rc\right)
\right\}\\
&-(1-\zeta h)\sigma\\
\geq{}& h\mathcal{F}\lc \widehat{c}_{j^*}(x^*),v_{j^*}(x^*)\rc +(1-\zeta h)\left(\lambda_{j^*}hv_{\bar \jmath^*}(x^*)
+(1-\lambda_{j^*}h)\mathbb{I}[v_{j^*}]\lc x^*+h\widehat{s}_{j^*}(x^*)\rc\right)
-(1-\zeta h)\sigma.
\end{aligned}
$$
For the last inequality above, we used $\widehat{c}_{j^*}(x^*)\in \mC^{\D}_{j^*}(x^*)$ which follows from \cref{Lemma: h=O(Dx) c}.
This yields
\beq\label{Eq: v-sigma}
\begin{aligned}
&\frac{\rho}{\theta}\lc v_{j^*}(x^*)-\sigma\rc\\
\geq{}& \mathcal{F}\lc \widehat{c}_{j^*}(x^*),v_{j^*}(x^*)\rc+\frac{\mathbb{I}[v_{j^*}]\lc x^*+h\widehat{s}_{j^*}(x^*)\rc-v_{j^*}(x^*)}{h}+\lambda_{j^*}(v_{\bar \jmath^*}(x^*)-v_{j^*}(x^*))+O(h)\\
\geq{}&\mathcal{F}\lc \widehat{c}_{j^*}(x^*),v_{j^*}(x^*)\rc+\widehat{s}_{j^*}(x^*)Dv_{j^*}(x^*)+\lambda_{j^*}(v_{\bar \jmath^*}(x^*)-v_{j^*}(x^*))+O(\D x).
\end{aligned}
\eeq
For the last inequality of \cref{Eq: v-sigma} we used \cref{Lemma: interp}. From \cref{Eq: v(x) sub} and \cref{Eq: v-sigma}, we then deduce $\zeta \sigma\leq O(\D x)$. \par
{\it{Step 2.}} Now let us define 
\beq\label{Eq: def sigma2}
\sigma_2=\max_{j,\, x\in \mathcal{G}^{\D x}} \{V^{\D}_j(x)-v_{j}(x)\}=V^{\D}_{j^*}(x^*)-v_{j^*}(x^*),\quad x^*\in \mathcal{G}^{\D x},
\eeq
and suppose $\sigma_2>0$. We know from \cref{Eq: def sigma2} that for all $j$ and $x\in \mathcal{G}^{\D x}$, $V^{\D}_j(x)\leq v_j(x)+\sigma_2$, hence
$$
V^{\D}_{\bar \jmath^*}(x^*)\leq v_{\bar \jmath^*}(x^*)+\sigma_2,\quad \mathbb{I}[V^{\D}_{j^*}]\lc x^*+h\widehat{s}_{j^*}(x^*)\rc\leq \mathbb{I}[v_{j^*}]\lc x^*+h\widehat{s}_{j^*}(x^*)\rc+\sigma_2,
$$
and we infer from the SL scheme \cref{Eq: HJBSL} that
$$
\begin{aligned}
&v_{j^*}(x^*)+\sigma_2 \\
\leq{}& \sup_{c\in \mC^{\D}_{j^*}(x^*)}\left\{
h\mathcal{F}\lc c,v_{j^*}(x^*)\rc +(1-\zeta h)\left(\lambda_{j^*}hv_{\bar \jmath^*}(x^*)
+(1-\lambda_{j^*}h)\mathbb{I}[v_{j^*}]\lc x^*+h(rx^*+y_{j^*}-c)\rc\right)\right\}\\
&+(1-\zeta h)\sigma_2.
\end{aligned}
$$
It follows from direct calculation that
\beq\label{Eq: v+sigma}
\begin{aligned}
&\frac{\rho}{\theta}\lc v_{j^*}(x^*)+\sigma_2\rc\\
\leq{}& \sup_{c\in \mC^{\D}_{j^*}(x^*)}\left\{\mathcal{F}\lc c,v_{j^*}(x^*)\rc+(1-\zeta h)(1-\lambda_{j^*}h)\frac{\mathbb{I}[v_{j^*}]\lc x^*+h(rx^*+y_{j^*}-c)\rc-v_{j^*}(x^*)}{h}\right\}\\
&+(1-\zeta h)\lambda_{j^*}(v_{\bar \jmath^*}(x^*)-v_{j^*}(x^*))\\
\leq {}&(1-\zeta h)(1-\lambda_{j^*}h)\sup_{c\in \mC^{\D}_{j^*}(x^*)}\left\{\mathcal{F}\lc c,v_{j^*}(x^*)\rc+\frac{\mathbb{I}[v_{j^*}]\lc x^*+h(rx^*+y_{j^*}-c)\rc-v_{j^*}(x^*)}{h}\right\}\\
{}&+\lambda_{j^*}(v_{\bar \jmath^*}(x^*)-v_{j^*}(x^*))+O(h).
\end{aligned}
\eeq
For the last inequality, we used the fact that $\mathcal{F}\lc c,v_{j^*}(x^*)\rc<0$ for all $c\in \mC^{\D}_{j^*}(x^*)$, hence $\mathcal{F}\lc c,v_{j^*}(x^*)\rc<(1-\zeta h)(1-\lambda_{j^*}h )\mathcal{F}\lc c,v_{j^*}(x^*)\rc.$\par
Next, we observe that since $v_{j^*}$ is concave, for all $c\in \mC^{\D}_{j^*}(x^*)$,
$$
\begin{aligned}
&\frac{\mathbb{I}[v_{j^*}]\lc x^*+h(rx^*+y_{j^*}-c)\rc-v_{j^*}(x^*)}{h}\\
\leq {}&\frac{v_{j^*}\lc x^*+h(rx^*+y_{j^*}-c)\rc-v_{j^*}(x^*)}{h}\leq (rx^*+y_{j^*}-c)Dv_{j^*}(x^*). 
\end{aligned}
$$
It follows then
\beq\label{Eq: F+interp}
\begin{aligned}
&\sup_{c\in \mC^{\D}_{j^*}(x^*)}\left\{\mathcal{F}\lc c,v_{j^*}(x^*)\rc+\frac{\mathbb{I}[v_{j^*}]\lc x^*+h(rx^*+y_{j^*}-c)\rc-v_{j^*}(x^*)}{h}\right\}\\
\leq {}&\sup_{c\in \mC^{\D}_{j^*}(x^*)}\left\{\mathcal{F}\lc c,v_{j^*}(x^*)\rc+(rx^*+y_{j^*}-c)Dv_{j^*}(x^*)\right\}\leq H\lc x^*,y_{j^*},v_{j^*}(x^*),Dv_{j^*}(x^*)\rc.
\end{aligned}
\eeq
We infer from \cref{Eq: v+sigma} and \cref{Eq: F+interp} that
$$
\begin{aligned}
&\frac{\rho}{\theta}\lc v_{j^*}(x^*)+\sigma_2\rc\\
\leq {}&(1-\zeta h)(1-\lambda_{j^*}h )H\lc x^*,y_{j^*},v_{j^*}(x^*),Dv_{j^*}(x^*)\rc+\lambda_{j^*}(v_{\bar \jmath^*}(x^*)-v_{j^*}(x^*))+O(h)\\
={}& H\lc x^*,y_{j^*},v_{j^*}(x^*),Dv_{j^*}(x^*)\rc+\lambda_{j^*}(v_{\bar \jmath^*}(x^*)-v_{j^*}(x^*))+O(h).
\end{aligned}
$$
For the last equality, we used the boundedness of $H\lc x^*,y_{j^*},v_{j^*}(x^*),Dv_{j^*}(x^*)\rc$. We observe that \cref{Eq: v(x) sub} holds here with $j^*,x^*$ defined by \cref{Eq: def sigma2}, yielding $\sigma_2\leq C\D x$. 
\end{proof}

\section{Solution method for an Aiyagari model with recursive utility}\label{Sec: MFG}
The wealth-income distribution of agents is given by a measure $m$ on $[\ux,\infty)\times \{y_1,y_2\}$ of the form $dm=\sum_{j\in\{1,2\}}dm_{j}(x)\otimes\delta_{y_j}(y)$, where $m_j$ denotes the distribution of wealth conditional on the
income state $y_j$. It is shown in \cite{achdou2026continuous} that $m_j$ has a density $g_j$ and possibly a Dirac mass at $x=\ux$ weighted by $\mu_j$. The pair of densities $(g_1,g_2)$ is described by the solution to a system of two stationary Fokker--Planck--Kolmogorov (FPK) equations in the sense of distributions:
\beq\label{Eq: FPK}
- \frac{d}{dx} \left[ s_j(x)g_j(x) \right] + \lambda_{\bar \jmath}g_{\bar \jmath}(x) - \lambda_j g_j(x) = 0,\quad \sum_j  \int_{x>\ux}g_j(x)dx+\sum_j  \mu_j=1.
\eeq
In \cref{Eq: FPK}, $s_j$ is the optimal saving policy, the transport term accounts for movements in wealth while the two switching terms account for transitions between income states. To solve \cref{Eq: FPK}, we follow \cite{SL_HA} and employ the dual SL scheme
\begin{equation}\label{Eq: FPK scheme}
G^{\D}_{i,j}
=
(1-\lambda_jh)
\sum_k
\beta_i\left(x_k+hs_{k,j}\right)G^{\D}_{k,j}
+
\lambda_{\bar\jmath}hG^{\D}_{i,\bar\jmath},
\qquad
\sum_{i}G^{\D}_{i,1}
+
\sum_{i}G^{\D}_{i,2}
=
1/{\D x},
\end{equation}
where $(\beta_i)_{i=0}^I$ denotes the piecewise-linear
$\mathbb{Q}_1$ nodal basis. This discretization depends only on the saving policy and the
income-state transition probabilities, and is therefore independent of
the particular discretization used for the HJB equation. Consequently,
the scheme \cref{Eq: FPK scheme} applies to both approaches
developed in this paper. We summarize the solution method for \cref{Eq: FPK scheme} in \cref{Algo: FPK}. 

\begin{algorithm}
\caption{Dual semi-Lagrangian scheme for the stationary FPK equation}
\label{Algo: FPK}
\begin{algorithmic}[1]
\Require optimal saving policies $s_{i,j}$, $j=1,2$, $i=0,\dots,I$
\State \textbf{Assemble transport matrices:} for $j=1,2$, build the sparse matrix
$M(s_j)\in\mathbb R^{(I+1)\times(I+1)}$
\State \textbf{Assemble the dual (transposed) operator:}
\[
Q=\begin{pmatrix}
(1-\lambda_1h)M^T(s_1)-I & \lambda_2h\,I\\[2pt]
\lambda_1h\,I & (1-\lambda_2h)M^T(s_2)-I
\end{pmatrix}
\]
\State \textbf{Normalization:} replace one row of $Q$ (e.g.\ the
first) by $\mathbf 1^T$ and set the corresponding right-hand side entry to
$1/\Delta x$; call the resulting system $\widetilde Q\,G=b$
\State \textbf{Solve} the sparse linear system $\widetilde Q\,G=b$ for
$G=(G_1,G_2)^T$
\State \textbf{Project:} $G\gets\max(G,0)$, then renormalize
$G\gets G\big/\big(\Delta x\sum_i(G_{i,1}+G_{i,2})\big)$
\State \Return $G_1,G_2$
\end{algorithmic}
\end{algorithm}
In \cref{Eq: FPK scheme}, the discrete density $G^{\D}_{i,j}$ approximates $g_j(x_i)$
for $x_i>\ux$, while $G^{\D}_{0,j}\D x$ approximates the mass $\mu_j$
at the borrowing constraint. For the implementation of the coupled HJB
and FPK equations, it is useful to express the interpolation operator in
matrix form. We define the matrix
$M(s_j)\in\mathbb{R}^{(I+1)\times(I+1)}$ by
$
M(s_j)_{i,k}
=
\beta_k\big(x_i+hs_{i,j}\big).
$
For a grid function $\varphi$, the matrix $M(s_j)$ represents the
composition of the interpolation operator with the shift
$x_i\mapsto x_i+hs_{i,j}$:
$
\big(M(s_j)\varphi\big)_i
=
\mathbb{I}[\varphi]\big(x_i+hs_{i,j}\big).
$
Its transpose satisfies
$
\big(M(s_j)^T\big)_{i,k}
=
\beta_i\big(x_k+hs_{k,j}\big),
$
so that
$
\big(M(s_j)^TG_j^\D\big)_i
=
\sum_k\beta_i\big(x_k+hs_{k,j}\big)G_{k,j}^\D.
$
Consequently, the transport term in \cref{Eq: FPK scheme} is represented
by the transpose $M(s_j)^T$.\par
For each $r$, solving the HJB equation yields the saving policy $s_j[r]$, which then gives the measure $m_j[r]$. We consider the aggregate wealth 
$
\mathcal{K}[r]=\sum_j\int_{x\geq\ux} xdm_j[r]. 
$
The Mean Field Game system for an Aiyagari model with the recursive utility \cref{EZ} is (cf. \cite[(1.7)]{achdou2026continuous}):
\begin{equation}\label{MFG}
\left\{\begin{aligned}
\quad &(i)\qquad &&0 = \max_{c\ge 0} \left\{ f(c,v_j) +(r^*x+y_j - c)Dv_j(x) \right\}+\lambda_j(v_{\bar \jmath}(x)-v_j(x)),\\
&&&c_{j}^*(x)=\mathop{\arg \max}\limits_{c\ge 0} \left\{ f(c,v_j) +(r^*x+y_j - c)Dv_j(x) \right\},\\
&(ii)&&- \frac{\partial}{\partial x} \left[ (r^*x+y_j - c_{j}^*(x)) g_j(x) \right]+\lambda_{\bar \jmath}g_{\bar \jmath}(x) -\lambda_jg_j(x)=0,\\
&&&\int_{x\geq \ux} g_1(x)dx+\int_{x\geq \ux} g_2(x)dx+\sum_{j\in \{1,2\}}\mu_j=1,\\
%&(iii_A)\qquad &&r^*=A\alpha \left(\frac{\mathcal{K}[r^*]}{N}\right)^{\alpha-1}-\delta=A\alpha \left(\frac{K[m^{(r^*)}]}{N}\right)^{\alpha-1}-\delta,\quad N=\dfrac{\lambda_2y_1+\lambda_1y_2}{\lambda_1+\lambda_2}.
\end{aligned}\right.
\end{equation}
with $r^*$ determined by the coupling condition
$$
(iii_A)\quad r^*=A\alpha\Big(\frac{\mathcal K[r^*]}{N}\Big)^{\alpha-1}-\delta=A\alpha \left(\frac{K[m^{(r^*)}]}{N}\right)^{\alpha-1}-\delta,\quad N=\frac{\lambda_2y_1+\lambda_1y_2}{\lambda_1+\lambda_2}.
$$
Here, $A$, $\alpha$, and $\delta$ denote, respectively, total factor productivity, the capital share in production, and the depreciation rate of capital. The equilibrium interest rate $r^*$ clears the capital market if $\mathcal{K}[r^*]=K_d(r^*)$, where $K(r)$ is the household capital supply obtained from the invariant distribution and $K_d(r)=N\big(\tfrac{A\alpha}{r+\delta}\big)^{1/(1-\alpha)}$ is the firm's capital demand, with aggregate labor $N$. \par
For each given $r$, we use an HJB equation solver to obtain an approximation of the saving policy $s_j(r)$, which gives us $G^{\D}_{i,j}$ from solving \cref{Eq: FPK scheme}. The supply of capital is approximated by 
$\mathcal{K}[r]=\sum_i x_i\big(G_{i,1}(r)+G_{i,2}(r)\big)\Delta x$.
 For solving \cref{MFG}, we define the excess capital supply $E(r) \;:=\; \mathcal{K}[r]- N\left(\frac{A\alpha}{r+\delta}\right)^{\frac1{1-\alpha}}$ and solve $E(r)=0$ by bisection. 

\section{Numerical examples}\label{Sec: numerical}
We consider the examples from \cite{achdou2026finite} and fix the following parameters:
\beq\label{Eq: parameters}
\begin{aligned}
&\rho=0.05,\,y_1=0.5,\,y_2=1.5,\,\lambda_1=0.2,\,\lambda_2=0.2,\,\ux=-0.15,\,\bar{x}=60,\\
&A=0.95,\,\alpha=0.35,\,\delta=0.1.
\end{aligned}
\eeq
The parameters in the first line of \cref{Eq: parameters} are used to solve the HJB equations, while those in the second line determine the equilibrium condition \cref{MFG} $(iii_A)$. Unless otherwise specified, we use the baseline discretization parameters $\Delta x=0.01$ and $h/\D x=0.2$. We consider the following three tests for the Aiyagari model (or, equivalently, the corresponding Mean Field Game), varying the risk-aversion coefficient $\gamma$ and the EIS coefficient $\psi$. The resulting equilibrium interest rates are reported below. For each test, we solve the HJB equation using \cref{Algo: Howard Newton} and solve the full model using the methods described in \cref{Sec: MFG}.
\begin{itemize}
\item[] {\textit{Test 1 MFG:}} $\g=2$, $\psi=0.4$, $r^*=0.0288$.
\item[] {\textit{Test 2 MFG:}} $\g=4$, $\psi=0.5$, $r^*=0.0266$.
\item[] {\textit{Test 3 MFG:}} $\g=20$, $\psi=0.5$, $r^*=0.0086$.
\end{itemize}
For comparison, the finite difference method in \cite{achdou2026finite} with the same parameters and grid size \(\Delta x=0.01\) yields equilibrium interest rates that agree with those reported above to four decimal places.\par
\begin{figure}[h!]
   \centering
    \caption{\cref{Algo: Howard Newton} sequences of grid functions in {\textit{Test 1 HJB}} (left) and {\textit{Test 2 HJB}} (right)}\label{Fig: vn_HN}
    	\begin{tabular}{cc}
	\includegraphics[width=0.46\textwidth]{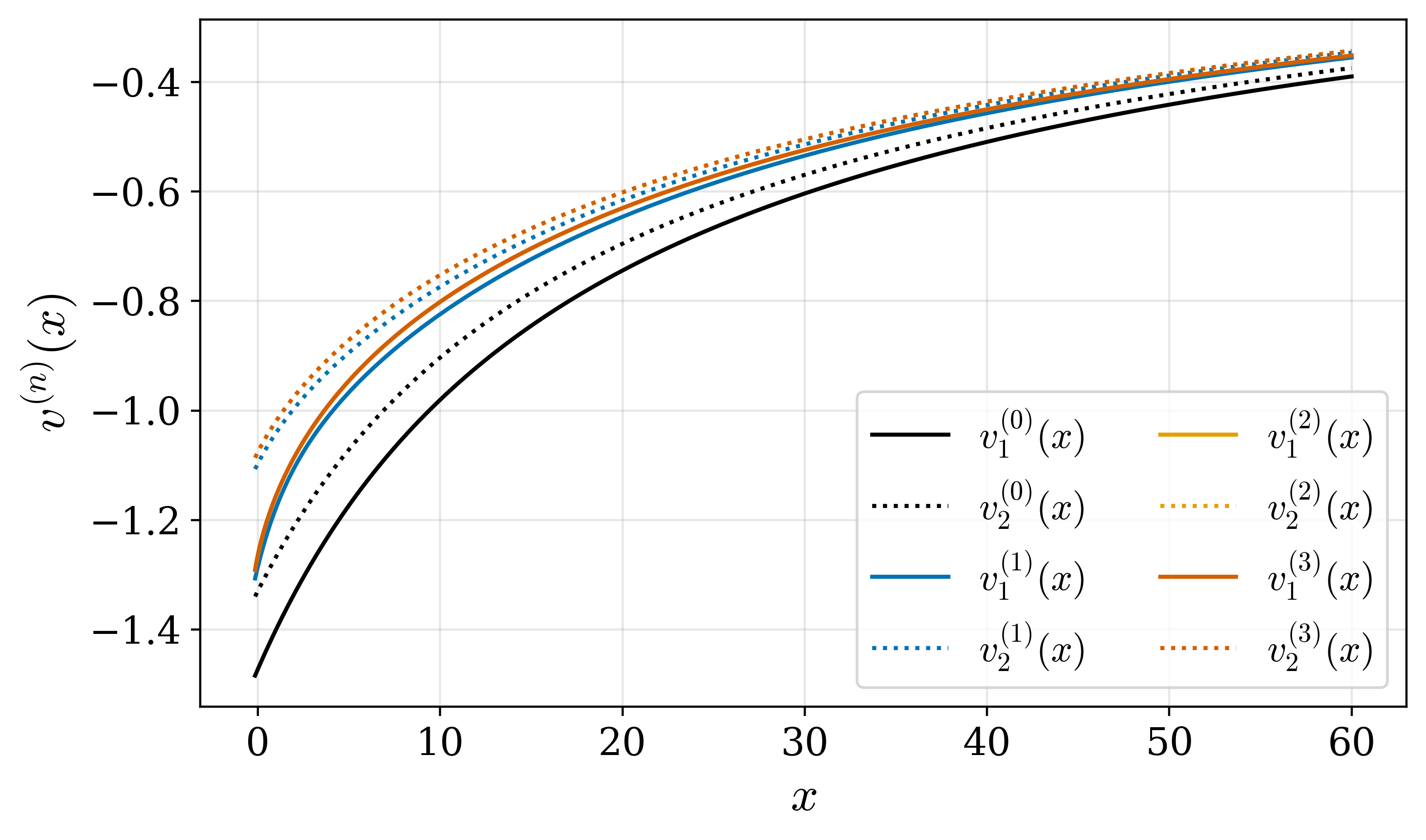} &
	\includegraphics[width=0.46\textwidth]{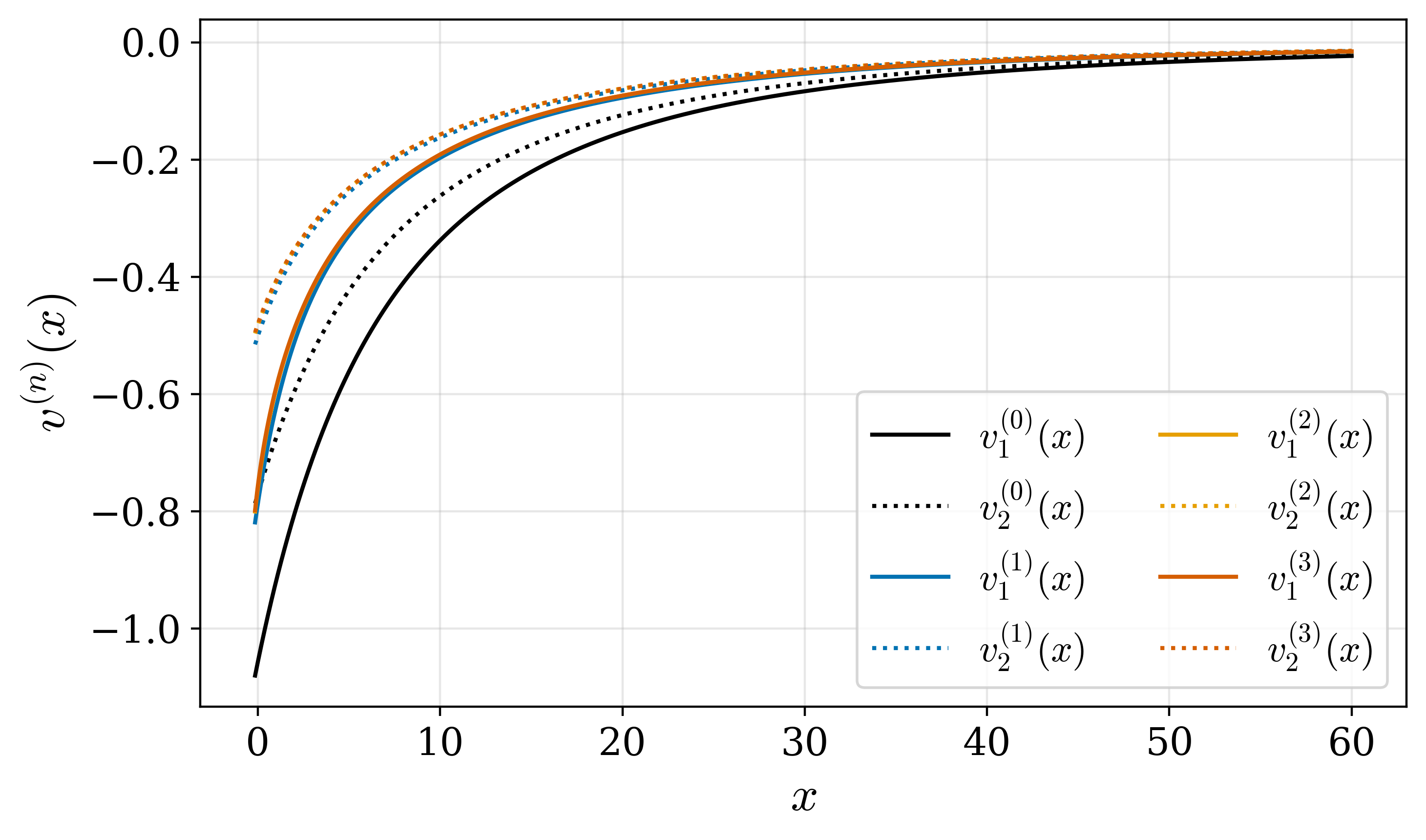} \\
	\end{tabular}
\end{figure}
\begin{figure}[h!]
   \centering
    \caption{\cref{Algo: up} constructs a non decreasing sequence of grid functions}\label{Fig: SL_VI}
    \includegraphics[width=0.5\linewidth]{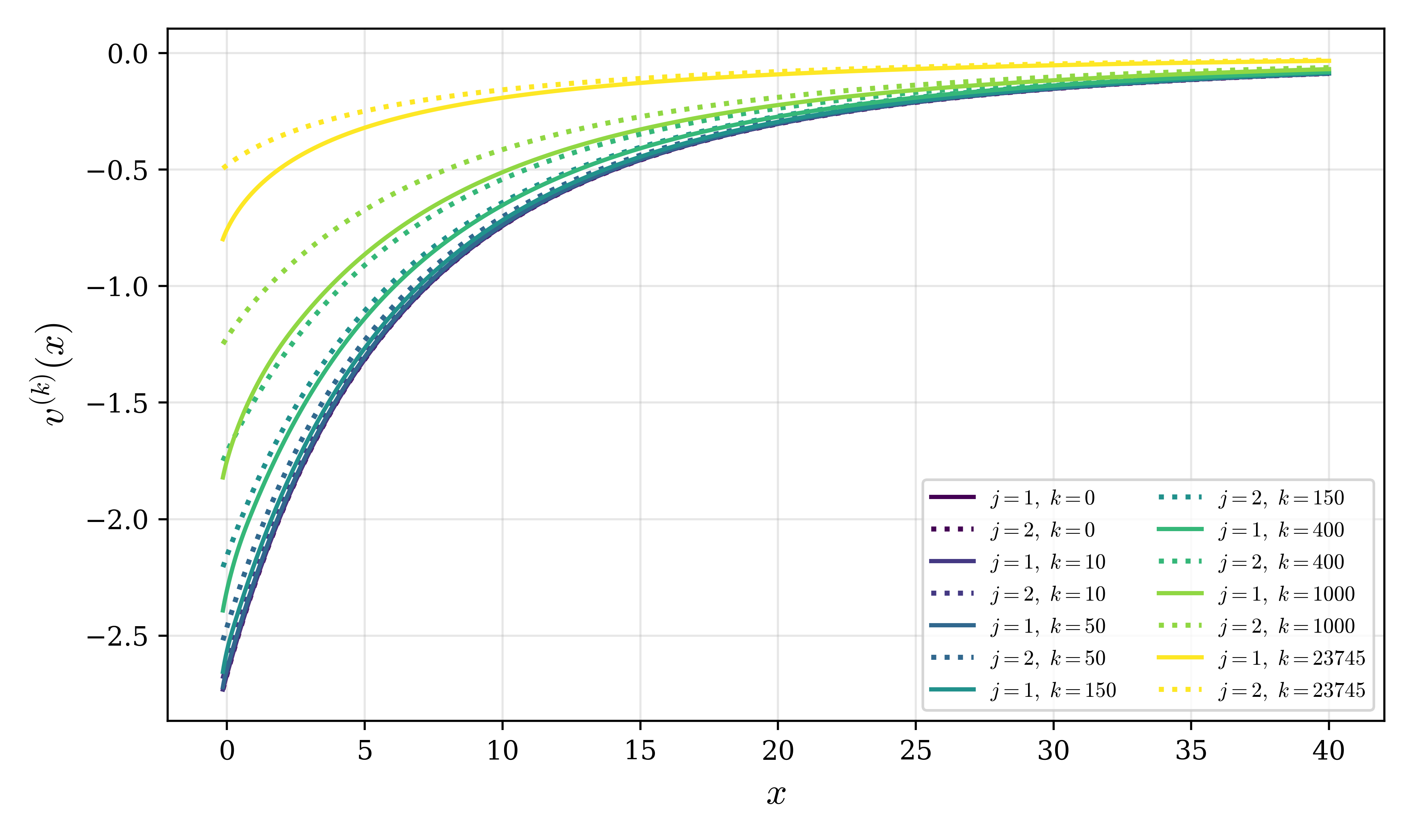}
\end{figure}

Next, we focus on the analysis of the HJB solver. To isolate these aspects from the equilibrium computation, we fix the interest rates at their corresponding equilibrium values. We refer to these  experiments as {\textit{Test 1 HJB}}, {\textit{Test 2 HJB}}, and {\textit{Test 3 HJB}}, respectively. \cref{Fig: vn_HN} illustrates that for {\textit{Test 1 and 2 HJB}}, the sequences produced by \cref{Algo: Howard Newton} are non decreasing as stated in \cref{Prop: convergence Howard}, i.e. $V^{(n)} \preceq V^{(n+1)}$. \cref{Fig: SL_VI} illustrates the fact, stated in \cref{Thm: SL_VI convergence}, that \cref{Algo: up} produces a non decreasing sequence of discrete functions for {\textit{Test 2 HJB}}. We also observe from \cref{Fig: vn_HN} and \cref{Fig: SL_VI} that both algorithms yield the same solution $V$ in {\textit{Test 2 HJB}}.

\begin{figure}[h!]
   \centering
    \caption{Convergence rate of \cref{Algo: Howard Newton}}\label{Fig: HN_rate}
    \includegraphics[width=0.4\linewidth]{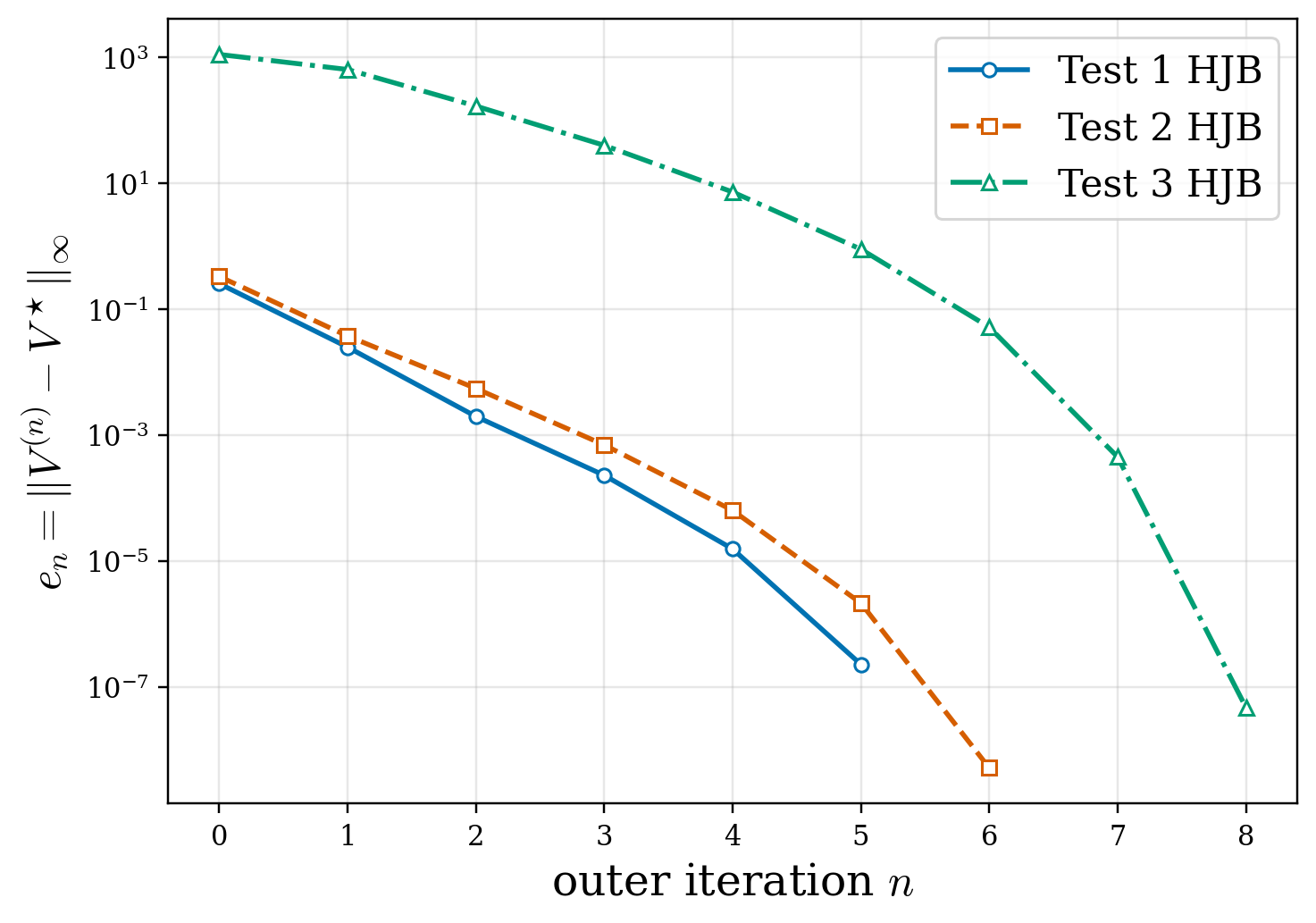}
\end{figure}

Figure~\ref{Fig: HN_rate} below plots the error $e_n=\|V^{(n)}-V^\star\|_\infty=\sup_{i,j}\vert V_{i,j}^{(n)}-V_{i,j}^\star \vert$ of the Howard--Newton iterates against the outer iteration $n$ (log scale), with $V^\star$ the last iterate. The convergence is superlinear, and numerically consistent with quadratic convergence in the last iterations. The quadratic convergence of the Howard algorithm for viscous HJB equations has been established in \cite{MR4392286,camilli2025quadratic}. The analysis of the convergence rate of \cref{Algo: Howard Newton} for the models of this paper, which have state constraints and use recursive utility, will be addressed in future work.  
\par
The numerical results show that the Howard--Newton algorithm converges significantly faster than value iteration (VI). The slow convergence of VI is mainly due to the small coefficient $\zeta h$. Although the SL scheme is not subject to the CFL condition and only requires $h=O(\Delta x)$ (cf. \cref{Thm: convergence rate}), the state constraint and accuracy considerations limit the practical choice of $h$. Compared with the baseline $h/\D x=0.2$, one may take $h/\Delta x=2$ when $\Delta x=0.001$ while obtaining comparable numerical results, but $h$ remains small. Therefore, the Howard--Newton algorithm is particularly advantageous for the models considered in this paper.\par

\begin{figure}[h!]
	\centering
	\caption{Convergence rate in {\textit{Test 1 HJB}} (left) and {\textit{Test 2 HJB}} (right)}\label{Fig: conv_SL_tests}
	\begin{tabular}{cc}
	\includegraphics[width=0.46\textwidth]{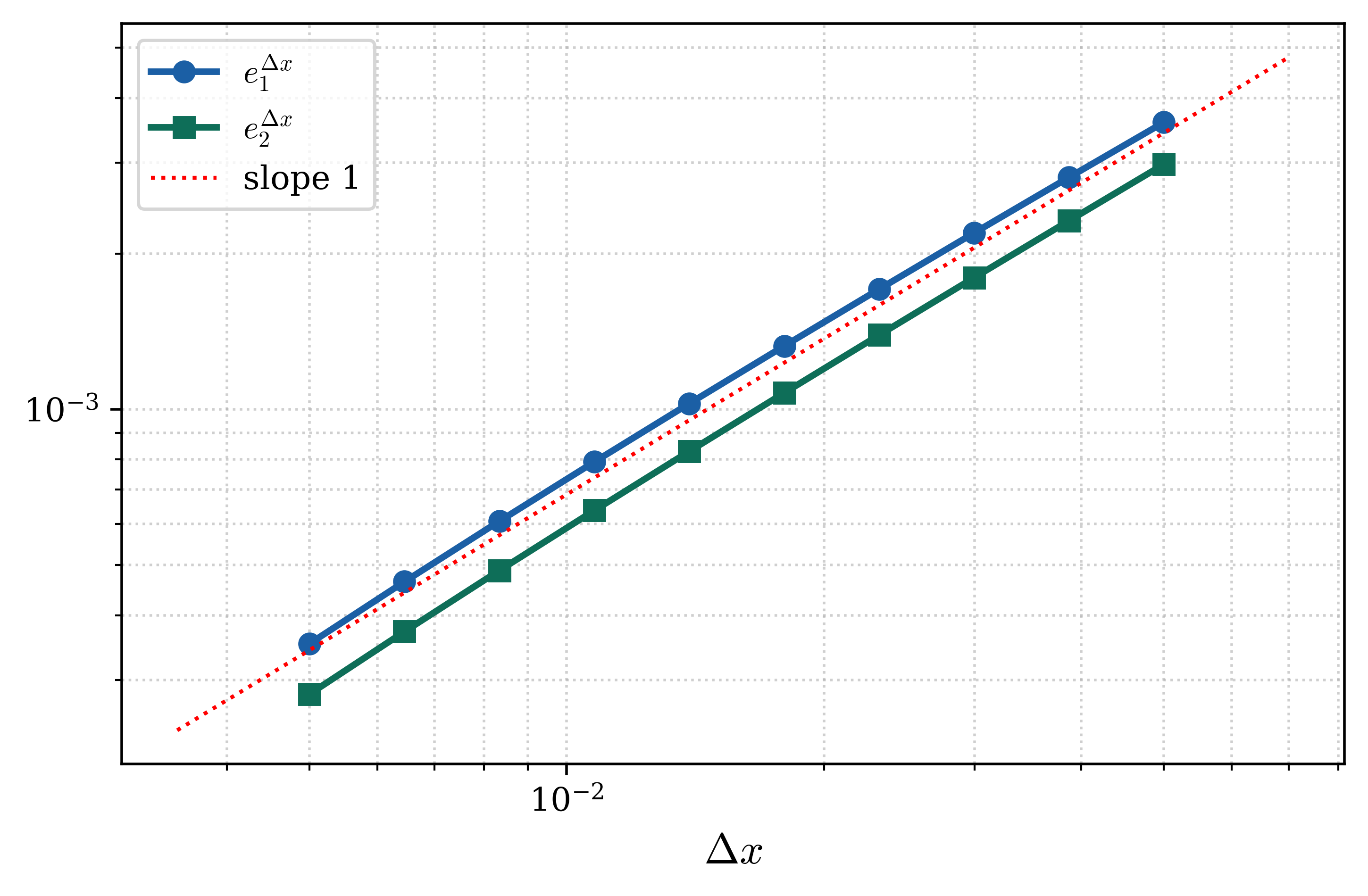} &
	\includegraphics[width=0.46\textwidth]{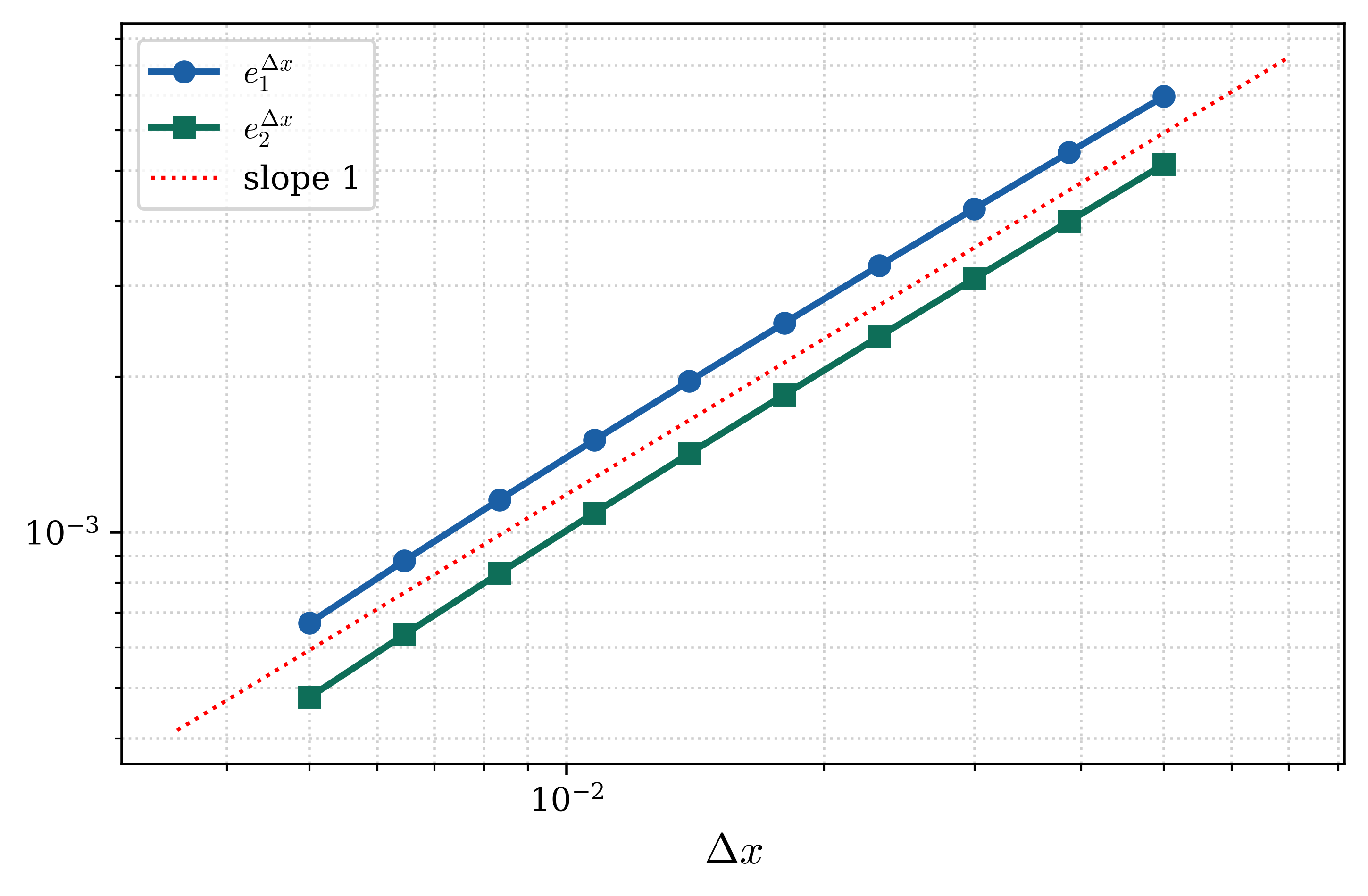} \\
	\end{tabular}
\end{figure}

We supplement a numerical estimation of the convergence rate. The results are obtained with \cref{Algo: Howard Newton}.  We compute a reference solution $(V_1^{\text{ref}}, V_2^{\text{ref}})$ on a fine grid with $\Delta x_{\text{ref}} = 0.0005$. We then solve the scheme on a sequence of ten coarser grids, with $\Delta x$ ranging from $0.005$ to $0.05$, spaced geometrically. For each coarse grid spacing $\Delta x > \Delta x_{\text{ref}}$, the reference solution is interpolated onto the coarse grid $\mathcal{G}^{\Delta x}$ via the piecewise linear
interpolation. We use $e_j^{\D x}=\max_{x\in \mathcal{G}^{\D x}} \left\vert \mathbb{I}[V_j^{\text{ref}}](x)-V_{j}(x)\right\vert$ to substitute $\max_{x\in \mathcal{G}^{\D x}} \left\vert v_j(x)-V_{j}(x)\right\vert$.  
\cref{Fig: conv_SL_tests} shows an approximately first-order convergence rate for both {\textit{Test 1 HJB}} (left) and {\textit{Test 2 HJB}} (right). For {\textit{Test 1 HJB}}, where $\theta>1$, the observed rate is consistent with the theoretical error bound established in \cref{Thm: convergence rate}. For {\textit{Test 2 HJB}}, where $0<\theta<1$, this provides numerical evidence that the same convergence-rate bound holds.

\begin{figure}[h!]
   \centering
    \caption{Asymptotic behavior of $s_1$ near $\ux$, in {\textit{Test 1 HJB}} (left) and {\textit{Test 2 HJB}} (right)}\label{Fig: s1}
    \includegraphics[width=0.8\linewidth]{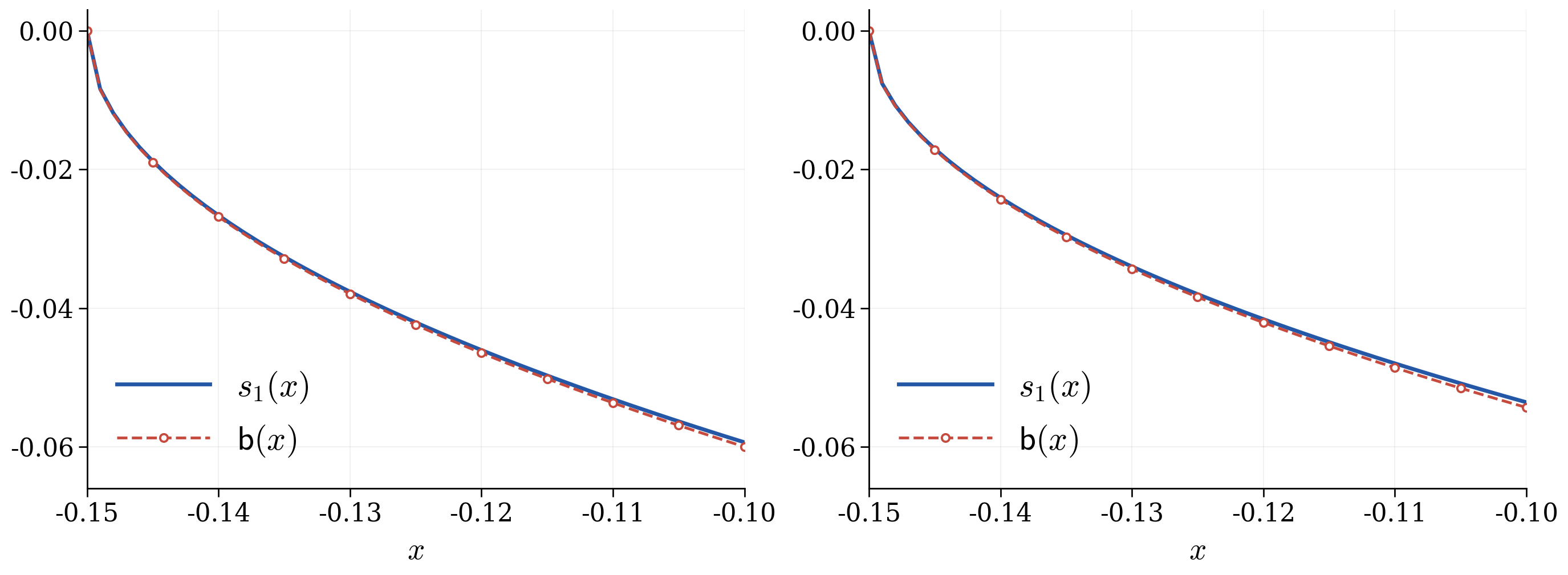}
\end{figure}

Finally, we numerically investigate the asymptotic behavior of \(s_1\) near the borrowing constraint \(\ux\), as described in \cref{Prop: Dv1 ux}. Let $\mathsf{b}(x)$ denote the leading-order term in \cref{Eq: -s1 bound}, namely,
$
\mathsf{b}(x)=-\sqrt{\frac{2\varkappa \p}{\rho} \lc r\ux+y_1\rc^{1+\p^{-1}}((1-\g)v_1(\ux))^{\frac{\p^{-1}-\g}{\g-1}}(x-\ux)},
$
where $\varkappa$ is given by \cref{Eq: def vark}. To evaluate $\varkappa$ numerically, we approximate $Dv_j(\ux)$ by the forward difference ${\mathbf{D}}^+(V_{0,j})$. In \cref{Fig: s1}, we compare the numerical saving policy $s_1(x)$ with the asymptotic approximation $\mathsf{b}(x)$ near $\ux$. For the numerical comparisons in \cref{Fig: s1}, we solve both tests with $h=\Delta x=0.001$ to improve the accuracy near the borrowing constraint.
\section{Concluding remarks}
An important direction of our ongoing research is to extend the semi-Lagrangian methods to continuous-time heterogeneous agent models with aggregate risk. In such models, the equilibrium can be characterized by a master equation, in which the value function depends not only on the individual state variables but also on the distribution of agents and the aggregate state. The resulting high-dimensional problem poses substantial analytical and computational challenges. Recent numerical approaches include finite difference methods combined with neural operators \cite{achdou2026mastering}, semi-Lagrangian methods combined with neural networks \cite{achdou2022simulating}, perturbation methods \cite{bilal2023solving}, global solution methods based on neural networks \cite{gu2024global} and reinforcement learning methods \cite{yang2025structural}. We plan to investigate how the methods developed in this paper can be extended to solving master equations, which may provide a framework for studying the interaction between aggregate shocks, wealth heterogeneity, and preferences for the timing of uncertainty resolution.
\appendix
\section{Accelerated exact solution of the policy-update problem}\label{App: argmax}
Evaluating \cref{Eq: c hat argmax k} requires a search over all indices $k$ at each grid point, at a cost of $O(I)$. This cost can be reduced by exploiting the concavity of the value function. The exact value function is concave, but the iterate $V^{(n)}_j$ need not be concave on the grid, although it typically is on large parts of it. We therefore partition the grid into regions on which $V^{(n)}_j$ is concave.\par
 \textbf{Concave runs.}
For $z\in [x_k,x_{k+1}]$, $\mathbb I[V_j](z)$ is affine in $z$ with slope $\mathbf D^+V_{k,j}$. A \emph{run} is a maximal set of consecutive indices $k$ for which $\mathbf D^+V_{k,j}$ is non-increasing in $k$; a new run starts at index $k$ whenever $\mathbf D^+V_{k,j}-\mathbf D^+V_{k-1,j}>\varepsilon$, with $\varepsilon=10^{-12}\max\{1,\max_{k}|\mathbf D^+V_{k,j}|\}$. This splits $\{0,\dots,I-1\}$ into runs $\mathcal R_1,\dots,\mathcal R_Q$, with $Q=1$ if $V_j$ is concave on the grid and $Q$ small if it fails to be concave only locally. On the union of the intervals $[x_k,x_{k+1}]$, $k\in\mathcal R_q$, $\Phi_{i,j}$ is concave in $c$ on the corresponding consumptions $\bigcup_{k\in\mathcal R_q}[\chi_{i,j}(x_{k+1}),\chi_{i,j}(x_k)]$.\par
 
 \textbf{Maximization on one run.}
On a run $\mathcal R_q$, let $k_q\in\mathcal R_q\cap K_{i,j}$ be an index such that
\beq\label{Eq: k q}
\max_{c\in\mC^{k_q}_j(x_i)}\Phi_{i,j}(c)=\max_{c\in\bigcup_{k\in\mathcal R_q}\mC^{k}_j(x_i)} \left\{\underbrace{h\mathcal F\lc c,V_{i,j}\rc +(1-\zeta h)(1-\lambda_jh)\,\mathbb{I}[V_j]\lc x_i+hs_{i,j}(c)\rc}_{\Phi_{i,j}(c)} \right\},
\eeq
that is, the consumptions $\mC^{k_q}_j(x_i)$ contain a maximizer of $\Phi_{i,j}$ over the run. Because $\Phi_{i,j}$ is concave in $c$ on the run, the sign of its derivative at the endpoints of the sets $\mC^{k}_j(x_i)$ changes monotonically with $k$. Hence $k_q$ can be found by bisection, using only the slopes $\mathbf D^+V_{k,j}$. The maximizer on the run is then $\widehat c^{\,k_q}_{i,j}$ from \cref{Eq: c hat k}, which is computed only once per run.\par
 
 \textbf{Global argmax.}
We compare $\Phi_{i,j}(\widehat c^{\,k_q}_{i,j})$ over $q=1,\dots,Q$ and take the largest (the smallest $q$ in case of ties). Since the runs cover all indices $k$, this gives the same $\widehat c_{i,j}$ as \cref{Eq: c hat argmax k}, up to the tolerance $\varepsilon$ and rounding. \Cref{Alg: Howard control update} summarizes the procedure.
\begin{algorithm}[h!]
\caption{Exact argmax for the consumption policy update \cref{Eq: Howard-Newton P-update SL}}
\label{Alg: Howard control update}
\begin{algorithmic}[1]
\Require $V_j=V^{(n)}_j$ on the grid, $j=1,2$; tolerance $\varepsilon$
\For{$j=1,2$}
  \State Split the indices $\{0,\dots,I-1\}$ into runs $\mathcal R_1,\dots,\mathcal R_Q$: a new run starts at index $k$ whenever $\mathbf D^+V_{k,j}-\mathbf D^+V_{k-1,j}>\varepsilon$
  \For{$i=0,\dots,I$}
    \For{$q=1,\dots,Q$}
      \State $k_q\gets$ an index in $\mathcal R_q\cap K_{i,j}$ satisfying \cref{Eq: k q}, found by bisection \Comment{concavity of $\Phi_{i,j}$ on $\mathcal R_q$}
      \State Compute $\widehat c^{\,k_q}_{i,j}$ from \cref{Eq: c hat k}
    \EndFor
    \State $q^\ast\gets$ smallest $q$ maximizing $\Phi_{i,j}\lc\widehat c^{\,k_q}_{i,j}\rc$; \quad $c^{(n+1)}_{i,j}\gets\widehat c^{\,k_{q^\ast}}_{i,j}$
  \EndFor
\EndFor
\end{algorithmic}
\end{algorithm}
\section{FOC-based approximate solution of the policy-update problem}\label{App: argmax FOC}
If the grid $\D=(h,\D x)$ is sufficiently fine, we expect the maximizer of the discrete problem to be well approximated by the control obtained from the continuous first-order condition \cref{c optimal}, with the derivative $Dv_j(x_i)$ replaced by ${\mathbf{D}}(V_{i,j}^{\Delta})$. This yields an approximate policy update method. \par 
\cref{Alg: FOC locator} uses the first-order-condition (FOC) as a locator to identify a candidate bracket $[x_{k_0},x_{k_0+1})$; if the raw gradient is nonpositive, no such guess is reliable, so we anchor the search at the household's own grid index instead. We then recompute consumption $\widehat{c}^{\,k}_{i,j}$ exactly on this bracket using its own secant slope as in \cref{Eq: c hat k}, and verify self-consistency by checking if indeed $x_k \le x_i + h\big(rx_i+y_j-\widehat{c}^{\,k}_{i,j}\big) \le x_{k+1}$ is satisfied. If it does not, we repeat the same closed-form calculation on a small number of neighboring brackets --- since in practice the correct bracket is almost always adjacent to the FOC's initial guess. Only when none of these candidates is self-consistent do we fall back to a brute-force search over the entire admissible grid, directly maximizing the right-hand side of the Bellman equation. The result may fail to be a global maximizer for \cref{Eq: c hat}, we therefore accept the candidate only where it does not decrease the one-step value.

\begin{algorithm}[h!]
\caption{FOC-based policy update with safeguard (approximate solution of \cref{Eq: Howard-Newton P-update SL})}
\label{Alg: FOC locator}
\begin{algorithmic}[1]
\Require $V_j=V^{(n)}_j$ and $c^{(n)}_j$ on the grid $x_0<\dots<x_I$, $j=1,2$; search radius $L$
\For{each grid point $(x_i, y_j)$}
    \If{$\mathbf{D}(V_{i,j}) > 0$}
        \State Set
        $$\tilde{c}_{i,j} \leftarrow \max\Bigg\{\min\left\{\rho^{\p}\big[\mathbf{D}(V_{i,j})\big]^{-\p}\big|(1-\gamma)V_{i,j}\big|^{\frac{1-\gamma \p}{1-\gamma}},\ \frac{x_i-\underline{x}}{h}+rx_i+y_j\right\},\ \frac{x_i-\bar x}{h}+rx_i+y_j\Bigg\}$$
        \State $k_0\gets$ the index with $x_{k_0} \le x_i + h(rx_i + y_j - \tilde c_{i,j}) < x_{k_0+1}$, clipped to $\{0,\dots,I-1\}$
    \Else
        \State $k_0 \gets \min\{i,I-1\}$ \Comment{anchor the search at the household's own grid index}
    \EndIf
    \State $k^\ast\gets\textsc{null}$; \quad $l\gets0$
    \While{$k^\ast=\textsc{null}$ \textbf{and} $l\le L$} \Comment{search cells in order of distance from $k_0$}
        \For{$k\in\{k_0-l,\ k_0+l\}$ with $0\le k\le I-1$ and $k^\ast=\textsc{null}$}
            \If{$\mathbf D^+(V_{k,j})>0$}
                \State Compute $\widehat{c}^{\,k}_{i,j}$ using \cref{Eq: c hat k}
                \If{$\widehat{c}^{\,k}_{i,j}>0$ \textbf{and} $x_k\le x_i+h\big(rx_i+y_j-\widehat{c}^{\,k}_{i,j}\big)\le x_{k+1}$}
                    \State $k^\ast\gets k$
                \EndIf
            \EndIf
        \EndFor
        \State $l\gets l+1$
    \EndWhile
    \If{$k^\ast \ne \textsc{null}$}
        \State $\widehat c_{i,j}\gets\widehat{c}^{\,k^\ast}_{i,j}$ \Comment{stationary point in cell $[x_{k^\ast},x_{k^\ast+1}]$}
    \Else \Comment{fallback: maximize over grid nodes}
        \State For each node $x_q$ with $\chi_{i,j}(x_q)=\frac{x_i-x_q}{h}+rx_i+y_j>0$, $q=0,\dots,I$, set
        $$
        q^\ast \gets \mathop{\arg\max}\limits_{q}\ \Phi_{i,j}\lc\chi_{i,j}(x_q)\rc,\quad \widehat c_{i,j}\gets \chi_{i,j}(x_{q^\ast})
        $$
    \EndIf
    \State \textbf{Safeguard:} $V^{+}_{i,j}\gets h\mathcal F\lc \widehat c_{i,j},V_{i,j}\rc +(1-\zeta h)\Bigl[\lambda_j hV_{i,\bjm}
  +(1-\lambda_jh)\mathbb{I}[V_j]\lc x_i+h(rx_i+y_j-\widehat c_{i,j})\rc \Bigr]$
    \If{$V^{+}_{i,j}\ge V_{i,j}$} $c^{(n+1)}_{i,j}\gets\widehat c_{i,j}$ \textbf{else} $c^{(n+1)}_{i,j}\gets c^{(n)}_{i,j}$ \EndIf
\EndFor
\end{algorithmic}
\end{algorithm}

\bigskip
\section*{Acknowledgments}
The author thanks Yves Achdou for helpful discussions and valuable suggestions for improving the paper.\par
Claude (Anthropic) and ChatGPT (OpenAI) assisted with some derivations, numerical checks of the formulas, and Python code that generates the quantitative illustrations. The author assumes responsibility for all content.
\small
\bibliographystyle{siamplain}
\bibliography{SL_recursive}

\end{document}